\documentclass[aos,MSNbibl,seceqn,nameyear,rotating,dvips]{arximspdf}
\usepackage{dcolumn}
\usepackage{graphicx}

\doi{10.1214/11-AOS952}
\volume{40}
\issue{1}
\pubyear{2012}
\firstpage{131}
\lastpage{158}

\makeatletter

\newcolumntype{d}[1]{D{.}{.}{#1}}

\newtheorem{theo}{Theorem}[section]

\newproclaim{rem}{Remark}
\newproclaim{condition}{Condition}

\renewcommand{\hat}{\widehat}

\newcommand{\ISE}{\mathrm{ISE}}
\newcommand{\var}{\operatorname{var}}
\newcommand{\as}{^*}
\newcommand{\bX}{{\bar X}}
\newcommand{\bbX}{{\mathbf X}}
\newcommand{\Bigmi}{\mi}

\newcommand{\cB}{{\mathcal B}}
\newcommand{\cI}{{\mathcal I}}
\newcommand{\cJ}{{\mathcal J}}
\newcommand{\cM}{{\mathcal M}}
\newcommand{\cP}{{\mathcal P}}
\newcommand{\cX}{{\mathcal X}}
\newcommand{\De}{\Delta}
\newcommand{\ep}{\varepsilon}
\newcommand{\half}{^{1/2}}
\newcommand{\hmu}{{\widehat\mu}}
\newcommand{\hp}{{\hat p}}
\newcommand{\hpi}{{\hat\pi}}
\newcommand{\inti}{\int_\cI}
\newcommand{\la}{\lambda}
\newcommand{\mhf}{^{-1/2}}
\newcommand{\mo}{^{-1}}
\newcommand{\mt}{^{-2}}
\newcommand{\mi}{\mid}
\newcommand{\ra}{\to}
\newcommand{\rai}{\ra\infty}
\newcommand{\si}{\sigma}
\newcommand{\sumi}{\sum_i}
\newcommand{\sumj}{\sum_j}
\newcommand{\thf}{{{1\over2}}}
\newcommand{\tmu}{{\tilde\mu}}
\newcommand{\tp}{{\tilde p}}

\makeatother

\begin{document}
\begin{frontmatter}

\title{Nonparametric regression with homogeneous group testing data\thanksref{T2}}
\runtitle{Nonparametric regression from grouped data}

\thankstext{T2}{Supported by grants and fellowships from the Australian Research Council.}

\begin{aug}
\author[A]{\fnms{Aurore} \snm{Delaigle}\corref{}\ead[label=e1]{A.Delaigle@ms.unimelb.edu.au}}
\and
\author[A]{\fnms{Peter} \snm{Hall}\ead[label=e2]{halpstat@ms.unimelb.edu.au}}
\runauthor{A. Delaigle and P. Hall}
\affiliation{University of Melbourne}
\address[A]{Department of Mathematics and Statistics\\
University of Melbourne, Parkville\\
Victoria, 3010\\
Australia\\
\printead{e1}\\
\phantom{E-mail: }\printead*{e2}} %adresu isvedimo komanda gale!
\end{aug}

% HISTORY:
\received{\smonth{4} \syear{2011}}
\revised{\smonth{10} \syear{2011}}

% ABSTRACT
%
\begin{abstract}
We introduce new nonparametric predictors for homogeneous pooled data
in the context of group testing for rare abnormalities and show that
they achieve optimal rates of convergence. In particular, when the
level of pooling is moderate, then despite the cost savings, the method
enjoys the same convergence rate as in the case of no pooling. In the
setting of ``over-pooling'' the convergence rate differs from that of
an optimal estimator by no more than a logarithmic factor. Our
approach improves on the random-pooling nonparametric predictor, which
is currently the only nonparametric method available, unless there is
no pooling, in which case the two approaches are identical.
\end{abstract}

% KEYWORDS
%
\begin{keyword}[class=AMS]
\kwd{62G08}.
\end{keyword}
\begin{keyword}
\kwd{Bandwidth}
\kwd{local polynomial estimator}
\kwd{pooling}
\kwd{prevalence}
\kwd{smoothing}.
\end{keyword}

\end{frontmatter}

%s1 ###
%s1 #&#
\section{Introduction}
In large screening studies where infection is detected by testing a
fluid (e.g., blood, urine, water, etc.), data are often pooled in groups
before the test is carried out, which permits savings in time and
money. This technique, known as group testing, dates back at least to
the Second World War, where \citet{Dor43} suggested using it to detect
syphilis in US soldiers. It has been used in a variety of large
screening studies, for example, to detect human immunodeficiency virus,
or HIV [\citet{GasHam89}], but pooling is also employed to
detect pollution, for example, in water or milk; see
\citet{NagRag}, \citet{Wahetal06}, \citet{Len07}, \citet{FahOurDeg06}. %
%Being able to do statistical inference from pooled data is
%particularly important for poor countries, where the lack of resources
%regularly prevents testing the blood of donors of relatively long
%periods. %See
%
%
Often in these studies, one or several explanatory variables are
available, in which case it is generally of interest to estimate the
conditional probability of infection. This problem has received
considerable attention in the group testing literature, where most
suggested techniques are parametric; see, for example,
\citet{VanGoeVer00}, \citet{BilTeb09} and
\citet{CheTebBil09}. Related
work includes that of \citet{CheSwa90}, \citet{GasJoh94},
\citet{HarPagSto98} and \citet{Xie01}.

Thus, although the original purpose of group testing was merely to
identify infected individuals more economically, the idea has since
been expanded extensively to include more general statistical
methodology when the data have to be gathered through grouping. Our
paper contributes in this context, developing and describing a
particularly effective approach to nonparametric regression. Obtaining
information in this way can be useful on its own, or for planning a
subsequent study.

Recently, \citet{DelMei11} suggested a nonparametric
estimator of the conditional probability of infection. Their method
enjoys optimal convergence rates when pooling is random, but it is not
consistent in the case of nonrandom, homogeneous pooling, which can be
defined as a setting where the covariates of individuals in a group
take similar values. In the parametric context it is well known that
homogeneous grouping improves the quality of estimators, but the
potential gains of homogeneous grouping are even greater in the
nonparametric context, where random grouping in moderate to large
groups can seriously degrade the quality of estimators.

We demonstrate that, when the data are grouped homogeneously, one can
construct more accurate nonparametric estimators of the conditional
probability of infection. We show that these improved estimators enjoy
faster, and optimal, convergence rates in a variety of contexts. Having
reliable estimators of the conditional probability of infection enables
more accurate identification of vulnerable categories of people, and
can lead to subsequent studies that can assist individuals who are
particularly vulnerable to infection. We illustrate the practical
performance of our procedure via simulated examples and an application
to the National Health and Nutrition Examination Survey (NHANES) study,
a large health and nutrition survey collected in the US; see
\href{http://www.cdc.gov/nchs/nhanes.htm}{www.cdc.gov/nchs/nhanes.htm}
for more about the NHANES research program.

%s2 ###
%s2 #&#
\section{Model and methodology}\label{sec2}

%s2.1 ###
%s2.1 #&#
\subsection{Main group testing model}

We observe independent and identically distributed (i.i.d.) data
$X_1,\ldots,X_N$, where $X$ is a covariate observed on each of $N$
respective objects (e.g., items or individuals), each of which is
subject to a potential, relatively rare ``abnormality.''
For example, $X$ could be the age or weight of an individual, and the
abnormality could be contamination by HIV.
Let $Y_i$ denote the result of a test on the $i$th object, such as
blood or urine test. That is, $Y_i$ takes the value 1 or 0 according to
whether the abnormality is detected or not, respectively. In large
screening studies, where $N$ is very large, testing each individual for
contamination can be too expensive or take too much time, and to
overcome this difficulty, it is common to pool data on several
individuals before performing the detection test.

Pooling is performed by partitioning the original dataset $\cX$,
comprised of the values $X_1,\ldots,X_N$, into $J$ subsets, or groups,
$\cX_1,\ldots,\cX_J$, say, where~$\cX_j$ is of size $n_j$ and
$n_1+\cdots+n_J=N$. We denote the elements of $\cX_j$ by
$X_{1j},\ldots,X_{n_jj}$. Each $X_{ij}$ corresponds to an $X_k$, and
each $X_k$ has a concomitant $Y_k$. If the $i$th element~$X_{ij}$ of
$\cX_j$ is $X_k$, then the concomitant of~$X_{ij}$ is $Y_{ij}=Y_k$.
Instead of trying to determine the value of $Y_{ij}$ directly, each
group~$\cX_j$ is tested to discover whether the abnormality is present
in the group, that is, to determine the value of
\[
Y_j\as=\max_{1\leq i\leq n_j} Y_{ij}
=\cases{
1, &\quad if $Y_{ij}=1$ for some $i$ in the range $1\leq i\leq n_j$,\cr
0, &\quad otherwise.}
\]
Of course,\vspace*{1pt} $Y_j\as$ is obtained without observing the
$Y_{ij}$'s directly; for example, when the abnormality is detected by a
blood test, the bloods of all individuals in a group are mixed
together, and this mixed blood is tested for contamination. From the
data pairs $(\cX _j,Y_j\as)$ we wish to estimate\vspace*{1pt} the probability
function $p(x)=P(Y_i=1\mi X_i=x)=E(Y_i=1\mi X_i=x)$.

Since $p$ is a regression curve, then if the sample $(X_i,Y_i),
i=1,\ldots,N$, were observed, we could use standard nonparametric
regression techniques such as, for example, local polynomial
estimators. Let $\ell\geq0$ be an integer, $h>0$ a bandwidth, $K$ a
kernel function and $K_h(x)=h^{-1}K(x/h)$. The standard $\ell$th
degree local polynomial estimator of $p$ is defined by
%
%e2.1 ###
%e2.1 #&#
\begin{equation}\label{LPp}
\hat p_{S}(x)=(1,0,\ldots,0)\mathbf{Q}\mo\mathbf{R},
\end{equation}
where
${\mathbf R}=(R_0(x),\ldots,R_\ell(x))^T$,
${\mathbf Q}=({\mathbf Q}_{ij})_{1\leq i,j \leq\ell+1}$, with ${\mathbf Q}_{ij}=Q_{i+j-2}(x)$, and where $Q_k(x)=\sum_{i=1}^N
(X_i-x)^kK_h(X_i-x)$ and $R_{k}(x)=\sum_{i=1}^N
Y_i(X_i-x)^kK_h(X_i-x)$. See, for example, \citet{FanGij96}.
Of course, when the data are pooled, the $Y_i$'s are not available, and
we cannot calculate such estimators. Therefore we need to develop
specific ways to estimate $p$ from pooled data.

%s2.2 ###
%s2.2 #&#
\subsection{Method for homogeneous pools}\label{sec21}
Depending on the study, it is not always possible to observe the $X_i$'s
before pooling the data, so that the individuals are pooled randomly.
This is the context of the work of \citet{DelMei11}, who
constructed a nonparametric estimator for the case where data $X_i$ are
assigned randomly to the groups $\cX_j$. See Appendix A.1 of the
supplemental article [\citet{DelHal}] for a summary of
properties of their estimator. In other studies, the $X_i$'s are
observed beforehand; see, for example, the study of hepatitis C
infection among 10,654 health care workers in Scotland, carried out by
Thorburn et~al. (\citeyear{Thoetal01}).
In such cases, it has already been demonstrated in the parametric
context that it can be greatly advantageous to pool the data
nonrandomly; see \citet{VanGoeVer00}.\vadjust{\goodbreak}

Unfortunately, the only nonparametric estimator available for group
testing data [see \citet{DelMei11}] crucially relies on
random grouping and is not valid when homogeneous groups are created.
Below we suggest a~new nonparametric approach which is valid with
homogeneous pooling. We introduce our procedure in the case of a single
covariate and equally sized groups. Generalizations of our method to
unequal group sizes and multiple covariates will be treated in Section
\ref{secgener}. These generalizations are similar in most respects.

To create homogeneous pools we divide the data into groups of equal
number, taking the $j$th group to be $\cX_j=\{X_{((j-1)\nu+1)},\ldots
,X_{(j\nu)}\}$, where $\nu=n_j$, in this case not depending on $j$,
is the number of data in each group, and $X_{(1)}\leq\cdots\leq
X_{(N)}$ denotes an ordering of the data in $\cX$. We assume that $\nu
$ divides $N$; the case where it does not is a particular case of our
generalization in Section \ref{secgener}.
Note that, with $Z_j\as=1-Y_j\as$,
%
%e2.2 ###
%e2.2 #&#
\begin{equation}\label{21}
E(Z_j\as\mi\cX)=\prod_{i=1}^\nu\{1-p(X_{ij})\} .
\end{equation}
The right-hand side here is generally close to $\{1-p(\bX_j)\}^\nu$,
where $\bX_j=\nu\mo\sumi X_{ij}$ denotes the average value of the
$X_{ij}$'s in the $j$th group, and that closeness motivates the
definition of $\hp(x)$ at (\ref{24}), below.
Let
%
%e2.3 ###
%e2.3 #&#
\begin{equation}\label{22}
\mu(x)=\{1-p(x)\}^\nu.
\end{equation}
Reflecting (\ref{21}) and the above discussion, we suggest
estimating $p(x)$ by
%
%e2.4 ###
%e2.4 #&#
\begin{equation}\label{24}
\hp(x)=1-\hmu(x)^{1/\nu} ,
\end{equation}
where $\hmu$ is a nonparametric estimator of $\mu$.

It remains to estimate $\mu$. We begin by giving motivation for our
methodology. Since, by construction, the groups are homogeneous, the
observations in a given group are similar. In particular,
$p(X_{((j-1)\nu+1)}),\ldots,p(X_{(j\nu)})$ are well approximated by
$p(\bX_j)$. Together, this and identity (\ref{21}) suggest that~$\mu(\bX_j)$
can be approximated by $E(Z_j\as\mi\bX_j)$, so that
$\mu(x)$ is approximately equal to the average of the $E(Z_j\as\mi
\bX_j)$'s over the $\bX_j$'s close to $x$, which can be estimated by
standard nonparametric regression estimators calculated from the data
$(\bX_j,Z_j\as)$, $j=1,\ldots,J$. Motivated by these considerations,
we define an $\ell$th order local polynomial estimator of $\mu$,
constructed from the data $(\bX_j,Z_j\as)$, by
%
%e2.5 ###
%e2.5 #&#
\begin{equation}\label{LPmu}
\hat\mu(x)=(1,0,\ldots,0){\mathbf S}\mo{\mathbf T},
\end{equation}
where ${\mathbf T}=(T_0(x),\ldots,T_\ell(x))^T$ and
${\mathbf S}=({\mathbf S}_{ij})_{1\leq i,j \leq\ell+1}$, with ${\mathbf S}_{ij}=S_{i+j-2}(x)$, $S_k(x)=\sumj(\bX_j-x)^kK_h(\bX_j-x)$, and
$T_{k}(x)=\sumj Z_j\as(\bX_j-x)^kK_h(\bX_j-x)$.\vspace*{1pt}

We shall show in Section \ref{secT} that this approach is well
founded, by proving consistency of the resulting estimator $\hat p$ of
$p$. We shall develop our theoretical results for a larger class of
estimators which encompasses the estimator at~(\ref{LPmu}).

%s3 ###
%s3 #&#
\section{Theoretical properties}\label{secT}
To study properties of our estimator it is convenient to express the
probability $p$, at a particular $x$, as
%
%e3.1 ###
%e3.1 #&#
\begin{equation}\label{T1}
p(x)=\delta(N)\pi(x) ,
\end{equation}
where $\delta=\delta(N)$ denotes a sequence of positive numbers that
potentially depend on~$N$, and $\pi$ is a fixed, nonnegative function.
To be as general as possible, we permit the group size $\nu=\nu
(N)\geq1$ to increase, and $\delta=\delta(N)$ to decrease, as $N$ diverges.

In large screening studies the abnormalities under investigation are
invariably rare, that is, $p$ is small. To understand the limitations
of our estimator, we shall study properties in the extreme situation
where $\delta\to0$ (and hence $p\to0$) as \mbox{$N\rai$}. More precisely,
we shall consider the ``low prevalence'' situation where
$
\nu\delta\ra0 %\label{T2}
$
as $N\rai$,
which is an asymptotic representation of the case where the group size
$\nu$ is relatively small and infection is rare.
In practice, groups larger than $10$ to $20$ are rarely taken. One
reason for this is that, depending on the proportion of positive
individuals in the population, some tests (e.g., HIV tests) become too
unreliable if the pool size is too large (larger than $\nu=5$ to 10 in
the HIV example). To reflect this fact, we shall also consider the
standard ``moderate pooling'' situation where
$\nu\delta\ra c>0 %\label{T19}
$
as $N\rai$.
However, there are tests for which groups could be taken as large as
$\nu=40$ to $50$. From the viewpoint of economics, large groups would
be beneficial, and might even be the only possible way to screen
individuals in poor countries. Hence we need to understand their
effects on the quality of estimators. We shall do this by investigating
asymptotic properties of our estimator in the extreme ``over-pooling''
situation where
$\nu\delta\rai$ %\label{T3}
as $N\rai$.

%In practice, groups are often kept relatively small so that estimators
%remain sufficiently accurate. We shall investigate the statistical
%properties of our estimator under two different pooling regimes. First
%we consider modest to moderate pooling, where the amount of pooling
%(i.e. the value of $\nu$) is chosen reasonably high, leading to cost
%savings when screening, but not so high that most of the groups $
%which case statistical accuracy deteriorates markedly. As we shall
%show, if $\nu\de$ is bounded as $N\rai$ then our estimator converges
%at rates equivalent to those in the case of no pooling. Much, but not
%all, of our asymptotic theory will treat the case where $\nu\de$ is
%bounded as $N$ increases, which of course encompasses settings where

%This corresponds to a modest amount of pooling, generally resulting
%from a relatively cautious approach.
%At the other extreme, estimators of $p$ computed from over-pooled
%data, where
%are generally relatively inaccurate (see section \ref{secT4}), and
%so over-pooling is to be avoided.

%s3.1 ###
%s3.1 #&#
\subsection{Conditions}

We shall derive theoretical properties of the estimator~$\hat p$
defined at (\ref{24}), where for $\hat\mu$ we shall generalize
the local polynomial estimators introduced at (\ref{LPmu}), by
considering a whole class of linear smoothers, defined by
%
%e3.2 ###
%e3.2 #&#
\begin{equation}\label{23}
\hmu(x)={\sumj w_j(x) Z_j\as\Big/\sumj w_j(x)} ,
\end{equation}
where\vspace*{1pt} the weights $w_j$ depend on $\cX$ but not on the
variables $Z_j\as$.
The local polynomial estimator defined at (\ref{LPmu}) can be
rewritten easily in this form, %; for example, using the notation
%$T_{k,j}(x)=(\bX_j-x)^kK_h(\bX_j-x)$, in the local constant case ($
and other popular nonparametric estimators (e.g., smoothing splines)
can be expressed in this form too; see, for example,
\citet{RupWanCar03}.

%In the local quadratic case we have
%$$w_j(x)=T_{0,j}S_2S_4-T_{0,j} S_3^2-T_{1,j}
%S_1S_4+T_{1,j}S_2S_3+T_{2,j}S_1S_3-T_{2,j}S_2^2$$
%so
%S_3+S_3^2S_1-S_3S_2^2=0
%and
%S_1S_4+S_2S_3^2+S_4S_1S_3-S_4S_2^2=0

Recall that $\bX_j=\nu\mo\sumi X_{ij}$ and let $h=h(N)$ denote a
sequence of constants decreasing to zero as $N\rai$. We can interpret
$h(N)$ as the bandwidth in a kernel-based construction of the weight
functions $w_j$ in (\ref{23}). Typically, the weights $w_j$ would
depend on $\bX_j$, and we assume that, for each $x\in\cI$, where
$\cI$ is a given compact, nondegenerate interval:
{\renewcommand{\thecondition}{S}
\begin{condition}\label{condS}

(S1) ${\sumj w_j(x) (\bX_j-x)/\sumj w_j(x)}=0$;

(S2) ${\sumj w_j(x) (\bX_j-x)^2/\sumj w_j(x)} =h^2
b(x)+o_p(h^2)$;

(S3) $\>{\sumj w_j(x)^2/\{\sumj w_j(x)\}^2}={\nu
v(x)/( Nh)}+o_p\{\nu/(Nh)\}$;

(S4) for\vspace*{1pt} each integer $k\geq1$, $\>{\sumj|w_j(x)|^k/\{
\sumj w_j(x)\}^k}=O_p[\{\nu/(Nh)\}^{k-1}]$,
where the functions $b$ and $v$ are continuous on $\cJ$ and are
related to the type of estimator.
\end{condition}}

%(a) the distribution of $X$ has a continuous density $f$ that is
%bounded away from zero on an open interval $\cJ$ containing $\cI$;
%(b) $p=\de\pi$ is bounded away from 1 uniformly in $x\in\cI$ and in $N
%vanish for $|\bX_j-x|>C h$, where $C>0$ is a constant; T.4
We also assume that:
{\renewcommand{\thecondition}{T}
\begin{condition}\label{condT}

{(T1)} the distribution of $X$ has a continuous density, $f$,
that is bounded away from zero on an open interval $\cJ$ containing
$\cI$;

{(T2)} $p=\delta\pi$ is bounded away from 1 uniformly in
$x\in\cI$ and in $N\geq1$;

{(T3)} the function $\pi$ in (\ref{T1}) has two H\"
older-continuous derivatives on $\cJ$;

{(T4)} for some $\ep>0$, $h+\nu\delta h+(\nu^2/N^{1-\ep
}h\delta)\ra0$ as $N\rai$;

{(T5)} the weights $w_j(x)$ vanish for $|\bX_j-x|>C h$,
where $C>0$ is a constant.
\end{condition}

The assumption in (T1) that $f$ is bounded away from zero on a compact
interval allows us to avoid pathological issues that arise when too few
values of $X$ are available in neighbourhoods of zeros of $f$. Finally,
when describing the size of $\hp(x)-p(x)$ simultaneously in many
values $x$ we shall ask that for some $C,\ep>0$,
%
%e3.3 ###
%e3.3 #&#
\begin{equation}\label{T6}\quad
\sup_{x, x'\in\cI\dvtx |x-x'|\leq N^{-C}} \biggl\{
{1\over|x-x'|^\ep}
\sum_k \biggl|{w_k(x)\over\sumj w_j(x)}-{w_k(x')\over\sumj
w_j(x')}\biggr|\biggr\}
=O_p(1) .
\end{equation}

For example, if the weights $w_j$ correspond to the local polynomial
estimator in (\ref{LPmu}) with $\ell=1$ (i.e., the local linear
estimator), with bandwidth~$h$ and a compactly supported, symmetric, H\"
older continuous, nonnegative kernel~$K$ satisfying \mbox{$\int K=1$}, and if
$h+(Nh)\mo=O(N^{-\ep_1})$ for some $\ep_1>0$, and (T1) holds, then
(T5), Condition \ref{condS} and~(\ref{T6}) hold with, in (S2) and (S3),
$b=\int u^2 K(u) \,du$ (not depending on~$x$) and $v(x)=f(x)\mo\int
K^2$. Furthermore, Condition \ref{condS} holds uniformly in $x\in\cI$. More generally it
is easy to see that when $\ell>1$, the $\ell$th order local
polynomial estimator in (\ref{LPmu}) satisfies $\sumj w_j(x) (\bX
_j-x)^k=0$ for $k=0,\ldots,\ell-1$, and hence conditions (S1) and~(S2)
are trivially satisfied. Conditions (S3) and (S4) too are
satisfied in this case, under mild conditions on the kernel. Note that
condition (S1) is not satisfied in the local constant case [$\ell=0$
in (\ref{LPmu})]. Although this instance can be easily accommodated by
modifying our conditions slightly, we simply omit it from our theory
because in practice the local linear estimator is almost invariably
preferred to the local constant one.
%
%re1 #&#
\begin{rem}\label{remlink}
Instead of linear smoothers, such as local polynomial estimators, we
could use alternative procedures which are sometimes preferred in the
context of binary dependent variables. For example,
\citet{FanHecWan95} suggest modeling the regression curve $m$ by
$m(x)=g^{-1}\{ \eta(x)\}$, where $g$ is a known link function, and
$\eta$ is an unknown curve. These methods have theoretical properties
similar to those of local polynomial estimators; the two methods differ
mostly through their bias, and, depending on the shapes of $m$ and $g$,
one method has a smaller bias than the other. We prefer local
polynomial estimators because they are easier to implement in practice.
\end{rem}

%s3.2 ###
%s3.2 #&#
\subsection{Low prevalence and moderate pooling}\label{seclow}
Our first result establishes convergence rates and asymptotic normality
for the estimator $\hat p$ defined at~(\ref{24}), with $\hat\mu$
at (\ref{23}). Note that we do not insist that $\nu$ and $\delta$
vary with~$N$; the regularity conditions for Theorem \ref{Theorem1}
hold in many cases where~$\nu$ and~$\delta$ are both fixed. Below we
use the notation $A(x)$ to denote the value taken by a function $A$ at
a point $x$, and the notation $A$ when referring to the function
itself. However, in some places, for example, in result (\ref{T7})
where it is necessary to refer explicitly to the point $x$ mentioned in
the statement ``for all $x\in\cI$,'' and in definitions (\ref
{T8}) and (\ref{T9}), where we are defining functions, the two
notations may appear a little ambiguous.
%
%th3.1 #&#
\begin{theo}\label{Theorem1}
Assume that Conditions \ref{condS} and \ref{condT} hold, and that $\nu\delta=O(1)$. Then,
for each $x\in\cI$,
%
%e3.4 ###
%e3.4 #&#
\begin{equation}\label{T7}
\hp(x)-p(x)=A(x) V(x)+B(x)+o_p\{\delta h^2+(\delta/Nh)\half
\} ,
\end{equation}
where the distribution of $V(x)$ converges to the standard normal law
as $N\to\infty$, and the functions $A$ and $B$ are given by
%
%e3.6 ###
%e3.5 ###
%e3.5 #&#
%e3.6 #&#
\begin{eqnarray}
\label{T8}
A&=&[(\nu Nh)\mo(1-p)^{2-\nu} \{1-(1-p)^\nu\}
v]\half=O\{(\delta/Nh)\half\}, \\
\label{T9}
B&=&\tfrac12 h^2 \{p''-(\nu-1) (1-p)\mo(p')^2\} b=O
(\delta h^2) ,
\end{eqnarray}
where $b$ and $v$ are as in \textup{(S2)} and \textup{(S3)}.
If, in addition, Condition \ref{condS} holds uniformly in $x\in\cI$, if (\ref
{T6}) holds, and if the functions $b$ and $v$ are bounded and
continuous, then
%
%e3.7 ###
%e3.7 #&#
\begin{equation}\label{T10}
\inti(\hp-p)^2=\inti(A^2+B^2)
+o_p\{\delta^2h^4+(\delta/Nh)\} .
\end{equation}
\end{theo}

Note that $A$ and $B$ represent, to first order, the standard deviation
of the error about the mean, and the main effect of bias, which arise
from the asymptotic distribution. For simplicity we shall call $A^2$
and $B$ the asymptotic variance and bias of the estimator.
%The remainder terms on the right-hand sides of \eqref{T7} and
%functions $\pi$, $b$ and $v$, in \eqref{T1}, (S2) and (S3), vanish.
%If those functions do not vanish at $x$ then \eqref{T7} can be
%replaced by $\hp(x)-p(x)=A(x) V(x)+\{1+o_p(1)\} B(x)$, and if
%neither function vanishes identically in $\cI$ then \eqref{T10} can
%be replaced by $\inti(\hp-p)^2=\{1+o_p(1)\}\inti(A^2+B^2)$.
From the theorem we see that, when $B(x)\neq0$ (e.g., for the local
polynomial estimator with $\ell=1$), if $N\delta\rai$ as $N\rai$,
then the rate of the estimator is optimized when $h$ is of size
$(N\delta)^{-1/5}$, in which case the estimator satisfies
%
%e3.8 ###
%e3.8 #&#
\begin{equation}\label{T17}
\mbox{for each $x\in\cI$} \qquad
\hp(x)-p(x)=O_p\{(\delta^3/N^2)^{1/5}\} .
\end{equation}
Note that when $\nu=1$ (no grouping), $\mu=1-p$ and our estimator of
$p$ reduces to a standard local linear smoother of $1-\mu$. For
example, the estimator at (\ref{LPmu}) coincides with $1-\hat p_{S}$
in (\ref{LPp}). Taking $\nu=1$ in the theorem, we deduce that the
convergence rate of our estimator for $\nu>1$, given at (\ref
{T17}), coincides with the rate for conventional linear smoothers
employed with nongrouped data.
By standard arguments it is straightforward to show that this rate is
optimal when $\pi$ has two derivatives, and hence our estimator is
rate optimal.
Although, in (T3), we assume that $\pi$ has two continuous
derivatives, continuity is imposed only so that the dominant term in an
expansion of bias can be identified relatively simply, and the
convergence rate at (\ref{T17}) can be derived without the
assumption of continuity. In addition, note that when $\nu\delta
=o(1)$ our estimator has the same asymptotic bias and variance
expressions, $B$ and $A$, as the estimator when $\nu=1$, which in that
case reduce to $A=(\delta/Nh)\half(\pi v)\half$ and $B=\thf
\delta h^2 \pi'' b+o_p(\delta h^2)$. In other words,
in that case the statistical cost of pooling is virtually zero.

The results discussed above also apply if performance is measured in
terms of integrated squared error (ISE), as at (\ref{T10}). In
particular, if $h$ is of size $(N\delta)^{-1/5}$, provided that $\nu
\delta$ is bounded, the estimator $\hp$ achieves the minimax optimal
convergence rate,
%
%e3.9 ###
%e3.9 #&#
\begin{equation}\label{T21}
\inti(\hp-p)^2=O_p\{(\delta^3/N^2)^{2/5}\} .
\end{equation}

%re2 #&#
\begin{rem}\label{remlocpoly}
Similar conclusions can be drawn in the case of estimators for which
$B(x)=0$, but this requires us to assume that the function $\pi$ has
enough derivatives so that an explicit, asymptotic, dominating, nonzero
bias term can be derived. For example, for our local polynomial
estimator of order $\ell>1$, we have $B(x)=0$ and the term $o_P(\delta
h^2)$ is only an upper bound to the bias of the estimator. A
nonvanishing asymptotic expression for the bias can easily be obtained
for $\ell>1$ if we assume that $\pi$ has $\ell+1$ continuous
derivatives. This can be done in a straightforward manner, but to keep
presentation simple, and since in practice local linear estimators are
almost invariably preferred to other local polynomial estimators, we
omit such expansions.
\end{rem}
%
%re3 #&#
\begin{rem}
In the case where $\delta\to0$, it could be argued that the rates are
meaningless since we are trying to estimate a function that tends to
zero, and that it is more appropriate to consider the nonzero part $\pi
$ of $p$ in the model at (\ref{T1}), and see how fast $\hpi=\hp
/\delta$ converges to $\pi$. The convergence rate of $\hat\pi$ is
easily deducible from (\ref{T17}):
%
%e3.10 ###
%e3.10 #&#
\begin{equation}\label{T18}
\mbox{For each $x\in\cI$} \qquad
\hpi(x)-\pi(x)=O_p\{(N\delta)^{-2/5}\} .\vadjust{\goodbreak}
\end{equation}
Provided that $N\delta\rai$ as $N\rai$, $\hpi(x)$ is consistent for
$\pi(x)$, and the convergence rate evinced by (\ref{T18}) is
optimal. % The estimator $\hpi$ is generally not of direct practical
%use, but if it is consistent then we can be assured that $\hp$ tracks
%the fluctuations of $p$, rather than simply capturing a measure of the
%smallness of that quantity.
\end{rem}

%s3.3 ###
%s3.3 #&#
\subsection{Over-pooling}\label{secOP}
The situation is quite different when $\nu\delta\rai$ as $N\rai$,
which can be interpreted as an asymptotic representation of the
situation where the data are pooled in groups of relatively large size
$\nu$. In practical terms the results in this section serve as a
salutary warning not to skimp on the testing budget. The work in
Section \ref{seclow} shows that the performance of estimators is robust, up to a
point, against increasing group size, but in the present section we
demonstrate that, after the dividing line between moderate pooling and
overpooling has been crossed, performance decreases sharply.

When $\nu\delta\rai$, properties of the estimator of $p(x)$ depend
on $x$, because there the order of magnitude of $\mu(x)$, at (\ref
{22}), depends critically on the rate at which $\{1-p(x)\}^\nu$
converges to zero.
%, or equivalently, $\{1-\de\pi(x)\}^\nu$, converges to zero; and
%this depends on the value of $\pi(x)$.
The following condition captures this aspect:
%
%e3.11 ###
%e3.11 #&#
\begin{equation}\label{T22}
\mbox{for some $\ep>0$} \qquad
\nu/h=o[N^{1-\ep} \{1-\delta\pi(x)\}^\nu] ,
\end{equation}
and the following theorem replaces Theorem \ref{Theorem1}.
%
%th3.2 #&#
\begin{theo}\label{Theorem2}
Assume that $\nu\delta\to\infty$ as $N\rai$, Conditions \ref{condS},
\ref{condT}\break
and~(\ref{T22}) hold, and $\pi$, $b$ and $v$ are all nonzero at $x$. Then
$
\hp(x)-p(x)=A(x) V(x)+\{1+o_p(1)\} B(x) ,
$
where $V(x)$ is asymptotically distributed as a~normal $N(0,1)$ as
$N\to\infty$, and $A$ and $B$ are given by the first identities in
each of (\ref{T8}) and (\ref{T9}).
\end{theo}

Note that the orders of magnitude given by the second identities in
each of (\ref{T8}) and (\ref{T9}) are not valid in this case,
and neither does result (\ref{T10}) necessarily hold under the
conditions of Theorem \ref{Theorem2}. %, since \eqref{T10} involves
%convergence at many different values of $x$.
Note too that the theorem can be extended to cases where $b=0$, along
the lines discussed in Remark \ref{remlocpoly}.
To elucidate the implications of Theorem \ref{Theorem2}, assume that~$\pi'(x)$
is nonzero, and define $\la_N(x)^5=\{1-\delta\pi(x)\}
^{-\nu}$, which, when $\nu\delta\to\infty$, diverges exponentially
fast as a function of $\nu\delta$. Given a sequence of constants
$c_N$ and a sequence of random variables $V_N$, write $V_N\asymp_p
c_N$ to indicate that both $V_N=O_p(c_N)$ and $c_N=O_p(V_N)$ as $N\rai
$. Theorem \ref{Theorem2} implies that, if $\nu\delta\to\infty$
and $h$ is a constant multiple of $\lambda_N (N\delta^4\nu
^3)^{-1/5}$, then
%
%e3.12 ###
%e3.12 #&#
\begin{equation}\label{T23b}
\{\hp(x)-p(x)\}^2\asymp_p(\delta^3/N^2)^{2/5} (\nu\delta)^{-2/5}
\la_N(x)^4,
\end{equation}
and in particular diverges at a rate that is exponentially slower, as a
function of~$\nu\delta$, than in the case where $\nu\delta=O(1)$,
treated in Section \ref{seclow}.
Result (\ref{T23b}) follows from the fact that $A(x)^2\asymp(\nu
Nh)\mo\la_N(x)^5$ and $|B(x)|\asymp h^2 \nu\delta^2$, where
$a_1(N)\asymp a_2(N)$ means that $a_1(N)/a_2(N)$ is bounded away from
zero and infinity. Note that (\ref{T23b}) includes the case where
$p$ (and hence $\delta$) is held fixed, and $\nu\rai$ as $N\rai$.\vadjust{\goodbreak}

The result at (\ref{T23b}) shows that when $\nu\rai$ as $N\rai$,
$\hp$ suffers from a~clear degradation of rates compared to the case
where $\nu\delta=O(1)$. Next we show that this degradation is
intrinsic to the problem, not to our estimator~$\hp$; any estimator
based on the pooled data in Section \ref{sec21} will experience an
exponentially rapid decline in performance as $\nu\delta\rai$.
More precisely we show in Theorem \ref{Theorem3} that, when $\nu
\delta\rai$ as $N\rai$, $\hp$ is near rate-optimal among all such
estimators.
Recall that, under our model (\ref{T1}), $p=\delta\pi$, where
$\delta=\delta(N)$ potentially converges to zero. If $\nu\delta\rai
$, then, by~(\ref{T23b}), we have\looseness=1
%
%e3.13 ###
%e3.13 #&#
\begin{equation}\label{T24}
|\hp(x)-p(x)|
=O_p[(\delta^3/N^2)^{1/5}
(\nu\delta)^{-1/5} \{1-p(x)\}^{-2\nu/5}] .
\end{equation}\looseness=0
Although this result was derived under the assumption that $\pi$ is a
fixed function with two continuous derivatives, since (\ref{T24})
is only an upper bound, then it is readily established under the
following more general assumption:\looseness=1
%
%e3.14 #&#
\begin{equation}
\label{T25}
\begin{tabular}{p{310pt}}
the nonnegative function $\pi=\pi_N$ can depend on $N$ and satisfies
$\pi_N(x) +$
$|\pi_N'(x)|+|\pi_N''(x)|\leq C_1$, for all $N$ and all $x$, where
the constant $C_1>0$
does not depend on $N$ or $x$.
\end{tabular}\hspace*{-28pt}
\end{equation}\looseness=0

Take the explanatory variables $X_i$ to be uniformly distributed on the
interval $\cM=[-\thf,\thf]$, and let $\cI\subset\cJ\subset\cM$
where 0 is an interior point of $\cI$. Let $p^1=\delta\pi_N$,
where $\pi_N$ satisfies (\ref{T25}), let $p^0\equiv\delta$
denote the version of $p^1$ when $\pi_N\equiv1$, and consider the condition
%
%e3.15 ###
%e3.15 #&#
\begin{equation}\label{T26}
(\nu^3\delta)^{1/2}=o\{N (1-\delta)^\nu\} .
\end{equation}
This assumption permits $\nu\delta$ to diverge with $N$, but not too
quickly. Indeed, using arguments similar to those in Section \ref
{secP3}, it can be shown that if (\ref{T26}) fails, then no
estimator of $p$ is consistent. Let $\cP$ be the class of measurable
functions ${\check p}$ of the pooled data pairs $(\cX_j,Y_j\as)$
introduced in Section \ref{sec21}.

%th3.3 #&#
\begin{theo}\label{Theorem3} Assume that $p^0$ and $p^1$ are bounded
below 1, that (\ref{T26}) holds and that $\nu\delta\rai$. Let
$x$ be an interior point of the support, $[-\thf,\thf]$, of the
uniformly distributed explanatory variables $X_i$. Then $C_2>0$, and
$\pi_N$, satisfying~(\ref{T25}), can be chosen such that
%
%e3.16 ###
%e3.16 #&#
\begin{eqnarray}\label{T27}
&&\liminf_{n\rai}
\max_{p=p^0, p^1} \inf_{{\check p}\in\cP} P[|{\check p}(x)-p(x)|\nonumber\\
&&\qquad\quad\hspace*{64pt}>C_2 \delta^{3/5}(N\nu\delta)^{-2/5}
\{1-p(x)\}^{-2\nu/5}]\\
&&\qquad>0 .\nonumber
\end{eqnarray}
\end{theo}

Except for the fact that $(\nu\delta)^{-2/5}$, rather than $(\nu
\delta)^{-1/5}$, appears in (\ref{T27}), the latter result
represents a converse to (\ref{T24}). The difference in powers here
is of minor importance since the main issue is the factor $\{1-p(x)\}
^{-2\nu/5}$, which (in the context $\nu\delta\rai$ of
over-pooling), diverges faster than any power of $\nu\delta$, and
this feature is represented in both (\ref{T24}) and
(\ref{T25}).

%s3.4 ###
%s3.4 #&#
\subsection{Comparison with the approach of Delaigle and Meister}

Arguments similar to those of \citet{DelMei11} can be used\vadjust{\goodbreak}
to show that, under conditions similar to those used in our Theorem
\ref{Theorem1}, their estimator~$\tp$ [see (A.1) in the supplemental
article, \citet{DelHal}] satisfies
$
\tp(x)-p(x)=A_1(x) V_1(x)+B_1(x)
+o_p\{\delta h^2+(\nu\delta/Nh)\half\} ,
$
where the random variable $V_1(x)$ has an asymptotic standard normal
distribution and\looseness=1
%
%e3.18 ###
%e3.17 ###
%e3.17 #&#
%e3.18 #&#
\begin{eqnarray}\quad
\label{T12}
A_1&=&[(Nh)\mo(1-p) q^{1-\nu}
\{1-(1-p) q^{\nu-1}\} v]\half
=O(\nu\delta/Nh) ,\\
\label{T13}
B_1&=&\tfrac12 h^2 p'' b
=O(\delta h^2)
\end{eqnarray}\looseness=0
with $q=E\{1-p(X)\}$. Likewise, the analog of (\ref{T10}) can be
derived in the following way:
$\inti(\tp-p)^2=\inti(A_1^2+B_1^2)
+o_p\{\delta^2h^4+(\nu\delta/Nh)\} .%\label{T14}
$
To simplify the comparison, assume that we use estimators for which $b$
and $v$ do not vanish, and that $\pi>0$.
We see when comparing (\ref{T12})--(\ref{T13}) with (\ref
{T8})--(\ref{T9}) that the asymptotic variance term $A^2$ of our
estimator is an order of magnitude~$\nu$ times smaller than $A_1^2$.
Note too the asymptotic bias terms of~$\hp$ and $\tp$ are of the same
size (the two biases are asymptotically equivalent if $\nu\delta\to
0$, and have the same magnitude in other cases). Hence, with our
procedure, the gain in accuracy can be quite substantial, especially if
$\nu$ is large.\looseness=-1

%s4 ###
%s4 #&#
\section{Numerical study}
We applied the local linear version of our local polynomial estimation
procedure [i.e., the one based on (\ref{LPmu}) with $\ell=1$] on
simulated and real examples. This method, which we denote below by DH,
is the one we prefer because it works well, it is very easy to
implement and we can easily derive and compute a good data-driven
bandwidth for it. The practical advantages of local linear estimators
over other local polynomial estimators have been discussed at length in
the standard nonparametric regression literature. %Among the many
%advantages is the simplicity of its implementation.
Of course, other versions of our general local linear smoother
procedure can be used, such as a spline approach or more complicated
iterative kernel procedures (see Remark \ref{remlink}). Each of the
methods gives essentially the same estimator.

In our simulations we compared the DH procedure, calculated by
definition from homogeneous groups, with the local linear estimator
$\hat p_{S}$ at (\ref{LPp})
that we would use if we had access to the original nongrouped data. We
also compared DH with the local linear version of the method of
\citet{DelMei11}, which, by definition, is calculated from
randomly created groups. We denote these two methods by LL and DM,
respectively. We took the kernel, $K$, equal to the standard normal
density. For $h$, in the DM case we used the plug-in bandwidth of
\citet{DelMei11} with their weight $\omega_0$; we used a
similar plug-in bandwidth in the LL and DH cases; see Section A.2 of
the supplemental article [\citet{DelHal}] for details.\looseness=1

%s4.1 ###
%s4.1 #&#
\subsection{Simulation results}\label{secsimres}
To facilitate the comparison with the DM method, we simulated data
according to the four models used by \citet{DelMei11}:

\begin{longlist}
\item
$p(x)=\{\sin(\pi x/2)+1.2\}/[ 20+40x^2\{\operatorname{sign}(x)
+1\}]$ and $X\sim U[-3,3]$ or $X\sim N(0,1.5^2)$;

\item $p(x)=\exp(-4+2x)/\{8+8\exp(-4+2x)\}$ and $X\sim U[-1,4]$ or
$X\sim N(2,1.5^2)$;

\item
$p(x)=x^2/8$ and $X\sim U[0,1]$ or $X\sim N(0.5,0.5^2)$;

\item
$p(x)=x^2/8$ and $X\sim U[-1,1]$ or $X\sim N(0,0.75^2)$.
\end{longlist}

We generated 200 samples from each model, with $X$ normal or uniform,
and with $N=1000$, $N=5000$ and $N=10\mbox{,}000$. Then for the DH method we
split each sample homogeneously into groups of equal sizes $\nu=5$,
$\nu=10$ or $\nu=20$; for the DM method, we created the groups
randomly (remember that this estimator is valid only for random groups).

\begin{sidewaystable}
\textwidth=\textheight
\tablewidth=\textwidth
\caption{Simulation results for models \textup{(i)} to \textup{(iv)},
when the $X_{i,j}$'s are uniform. The numbers show \mbox{$10^4$ $\times$ MED}
(IQR) of the ISE calculated from 200~simulated~samples}
\label{tableunif}
\begin{tabular*}{\tablewidth}{@{\extracolsep{\fill}}lcd{1.10}d{2.10}d{2.10}d{2.10}d{2.10}d{2.10}c@{}}
\hline
& & \multicolumn{1}{c}{\hspace*{-3pt}$\bolds{\nu=1}$} & \multicolumn{2}{c}{\hspace*{-3pt}$\bolds{\nu=5}$}
& \multicolumn{2}{c}{\hspace*{-3pt}$\bolds{\nu=10}$}
& \multicolumn{2}{c@{}}{$\bolds{\nu=20}$} \\[-4pt]
& & \multicolumn{1}{l}{\hrulefill\hspace*{3pt}}
& \multicolumn{2}{l}{\hrulefill\hspace*{3pt}}
& \multicolumn{2}{l}{\hrulefill\hspace*{3pt}}
& \multicolumn{2}{l@{}}{\hrulefill}\\
\textbf{Model} & $\bolds{N}$ & \multicolumn{1}{c}{\hspace*{-3pt}\textbf{LL}} & \multicolumn{1}{c}{\hspace*{-3pt}\textbf{DH}}
& \multicolumn{1}{c}{\hspace*{-3pt}\textbf{DM}}
& \multicolumn{1}{c}{\hspace*{-3pt}\textbf{DH}}
& \multicolumn{1}{c}{\hspace*{-3pt}\textbf{DM}} & \multicolumn{1}{c}{\textbf{DH}}
& \multicolumn{1}{c@{}}{\textbf{DM}} \\
\hline
\hphantom{ii}(i)
&$10^3$&
9.35\mbox{ (7.42)}&10.1\mbox{ (8.16)}&26.9\mbox{ (24.1)}&11.0\mbox{ (8.53)}&51.2\mbox{ (49.6)}&17.8\mbox{ (484)}& 122\mbox{ (110)}\hphantom{.00}\\
&$5\cdot10^3$
&2.91\mbox{ (2.01)}&2.94\mbox{ (2.38)}&7.59\mbox{ (5.34)}&3.30\mbox{ (2.06)}&14.1\mbox{ (11.4)}&4.46\mbox{ (2.94)}&\hphantom{0}29.2\mbox{ (25.2)}\hphantom{0}\\
&$10^4$&
1.62\mbox{ (1.20)}&1.83\mbox{ (1.40)}&4.54\mbox{ (3.05)}&2.07\mbox{ (1.63)}&7.70\mbox{ (6.13)}&2.89\mbox{ (1.95)}&\hphantom{0}16.8\mbox{
(13.9)}\hphantom{0}
\\[3pt]
\hphantom{i}(ii)
&$10^3$&
6.37\mbox{ (8.38)}&8.66\mbox{ (9.99)}&29.4\mbox{ (28.4)}&10.3\mbox{ (11.4)}&64.7\mbox{ (69.5)}&29.7\mbox{ (1560)}& 166\mbox{
(169)}\hphantom{.00}
\\
&$5\cdot10^3$&
1.48\mbox{ (1.37)}&1.66\mbox{ (2.26)}&6.37\mbox{ (5.93)}&2.41\mbox{ (2.74)}&13.8\mbox{ (12.1)}&4.47\mbox{ (5.94)}&\hphantom{0}35.8\mbox{ (30.0)}\hphantom{0}\\
&$10^4$&
0.963\mbox{ (0.843)}&1.02\mbox{ (1.16)}&3.39\mbox{ (2.89)}&1.35\mbox{ (1.25)}&7.04\mbox{ (6.20)}&2.35\mbox{ (3.26)}&\hphantom{0}19.1\mbox{ (17.2)}\hphantom{0}\\
[3pt]
(iii)
&$10^3$&
0.777\mbox{ (0.978)}&0.860\mbox{ (1.26)}&3.44\mbox{ (4.03)}&1.02\mbox{ (1.31)}&7.26\mbox{ (8.37)}&1.90\mbox{ (4.81)}&\hphantom{0}19.9\mbox{
(19.5)}\hphantom{0}
\\
&$5\cdot10^3$&
0.176\mbox{ (0.220)}&0.166\mbox{ (0.254)}&0.722\mbox{ (0.818)}&0.214\mbox{ (0.298)}&1.68\mbox{ (1.67)}&0.356\mbox{ (0.482)}&\hphantom{00}4.48\mbox{ (3.97)}\\
&$10^4$&
0.093\mbox{ (0.108)}&0.100\mbox{ (0.128)}&0.355\mbox{ (0.344)}&0.117\mbox{ (0.158)}&0.797\mbox{ (0.800)}&0.200\mbox{ (0.212)}&\hphantom{00}2.28\mbox{ (1.79)}
\\[3pt]
\hspace*{1pt}(iv)
&$10^3$&
2.33\mbox{ (2.11)}&2.49\mbox{ (2.32)}&7.41\mbox{ (9.81)}&2.70\mbox{ (2.55)}&17.2\mbox{ (16.3)}&5.07\mbox{ (166)}&\hphantom{0}39.7\mbox{ (34.1)}\hphantom{0}\\
&$5\cdot10^3$&
0.590\mbox{ (0.510)}&0.633\mbox{ (0.602)}&2.01\mbox{ (1.73)}&0.637\mbox{ (0.702)}&4.05\mbox{ (3.70)}&0.964\mbox{ (1.06)}&\hphantom{00}9.62\mbox{ (9.11)}\\
&$10^4$&
0.309\mbox{ (0.254)}&0.317\mbox{ (0.293)}&1.10\mbox{ (0.873)}&0.373\mbox{ (0.311)}&2.31\mbox{ (1.89)}&0.570\mbox{ (0.539)}&\hphantom{00}5.47\mbox{ (4.80)}\\
\hline
\end{tabular*}
\end{sidewaystable}

To assess the performance of our DH estimator we calculated, in each
case and for each of the 200 generated samples, the integrated squared
error $\ISE=\int_a^b (\hat p-p)^2$, with $a$ and $b$ denoting the
0.05 and 0.95 quantiles of the distribution of $X$. We did the same for
the DM and LL estimators $\tilde p$ and~$\hat p_S$. For brevity,
figures illustrating the results are provided in Section A.4 of the
supplemental article [\citet{DelHal}], and here we show only
summary statistics. In the graphs of Section A.4, we show the target
curve (thin uninterrupted curve) as well as three interrupted curves;
these were calculated from the samples that gave the first, second and
third quartiles of the 200 ISE values.

%t2 ###
%t2 #&#
\begin{sidewaystable}
\textwidth=\textheight
\tablewidth=\textwidth
\caption{Simulation results for models \textup{(i)} to \textup{(iv)}, when the
$X_{i,j}$'s are normal. The numbers show $10^4$ $\times$ MED (IQR) of the
ISE calculated from 200~simulated~samples}
\label{tablenorm}
\begin{tabular*}{\tablewidth}{@{\extracolsep{\fill}}lc
d{2.10}d{2.10}d{2.10}d{2.10}d{2.8}d{2.10}c@{}}
\hline
& & \multicolumn{1}{c}{\hspace*{-3pt}$\bolds{\nu=1}$} & \multicolumn{2}{c}{\hspace*{-3pt}$\bolds{\nu=5}$}
& \multicolumn{2}{c}{\hspace*{-3pt}$\bolds{\nu=10}$}
& \multicolumn{2}{c@{}}{$\bolds{\nu=20}$} \\[-4pt]
& & \multicolumn{1}{l}{\hrulefill\hspace*{3pt}}
& \multicolumn{2}{l}{\hrulefill\hspace*{3pt}}
& \multicolumn{2}{l}{\hrulefill\hspace*{3pt}}
& \multicolumn{2}{l@{}}{\hrulefill}\\
\textbf{Model} & $\bolds{N}$ & \multicolumn{1}{c}{\hspace*{-3pt}\textbf{LL}} & \multicolumn{1}{c}{\hspace*{-3pt}\textbf{DH}}
& \multicolumn{1}{c}{\hspace*{-3pt}\textbf{DM}}
& \multicolumn{1}{c}{\hspace*{-3pt}\textbf{DH}}
& \multicolumn{1}{c}{\hspace*{-3pt}\textbf{DM}} & \multicolumn{1}{c}{\textbf{DH}}
& \multicolumn{1}{c@{}}{\textbf{DM}} \\
\hline
\hphantom{ii}(i)
&$10^3$&
10.3\mbox{ (6.69)}&10.7\mbox{ (7.18)}&20.8\mbox{ (19.0)}&10.8\mbox{ (8.04)}&37.0\mbox{ (35.3)}&12.8\mbox{ (9.70)}&\hphantom{0}85.6\mbox{ (72.8)}\hphantom{0}\\
&$5\cdot10^3$&
4.35\mbox{ (2.80)}&4.14\mbox{ (2.71)}&9.60\mbox{ (5.49)}&4.32\mbox{ (2.95)}&12.0\mbox{ (11.1)}&4.50\mbox{ (3.44)}&\hphantom{0}17.3\mbox{ (18.8)}\hphantom{0}\\
&$10^4$&
3.12\mbox{ (1.77)}&3.33\mbox{ (2.07)}&7.66\mbox{ (4.12)}&3.01\mbox{ (2.01)}&9.42\mbox{ (5.68)}&3.20\mbox{ (2.19)}&\hphantom{0}13.6\mbox{ (11.0)}\hphantom{0}\\
\hphantom{i}(ii)
&$10^3$&
5.02\mbox{ (5.20)}&5.78\mbox{ (6.83)}&17.0\mbox{ (23.0)}&8.18\mbox{ (10.6)}&46.0\mbox{ (57.8)}&21.1\mbox{ (64.0)}&167
(202)\hphantom{00.}
\\
&$5\cdot10^3$&
1.69\mbox{ (1.95)}&1.98\mbox{ (2.18)}&4.23\mbox{ (5.97)}&2.36\mbox{ (3.40)}&9.48\mbox{ (12.3)}&5.37\mbox{ (6.75)}&\hphantom{0}28.3\mbox{
(36.9)}\hphantom{0}
\\
&$10^4$&
1.02\mbox{ (0.925)}&1.17\mbox{ (1.21)}&2.99\mbox{ (3.12)}&1.46\mbox{ (1.64)}&5.51\mbox{ (6.81)}&3.04\mbox{ (3.22)}&\hphantom{0}15.0\mbox{
(17.7)}\hphantom{0}
\\
(iii)
&$10^3$&
0.897\mbox{ (1.53)}&0.885\mbox{ (1.06)}&2.95\mbox{ (3.36)}&0.910\mbox{ (1.27)}&5.73\mbox{ (7.10)}&1.37\mbox{ (2.14)}&\hphantom{0}23.7\mbox{
(27.3)}\hphantom{0}
\\
&$5\cdot10^3$&
0.274\mbox{ (0.389)}&0.263\mbox{ (0.325)}&0.946\mbox{ (0.997)}&0.260\mbox{ (0.383)}&1.61\mbox{ (2.08)}&0.448\mbox{ (0.692)}&\hphantom{00}4.26\mbox{
(4.93)}
\\
&$10^4$&
0.204\mbox{ (0.270)}&0.148\mbox{ (0.175)}&0.637\mbox{ (0.725)}&0.182\mbox{ (0.219)}&1.13\mbox{ (1.10)}&0.323\mbox{ (0.435)}&\hphantom{00}2.42\mbox{
(2.58)}
\\
\hspace*{1pt}(iv)
&$10^3$&
4.13\mbox{ (4.30)}&3.60\mbox{ (3.48)}&13.2\mbox{ (12.5)}&4.32\mbox{ (3.84)}&28.1\mbox{ (26.9)}&7.60\mbox{ (9.43)}&\hphantom{0}82.3\mbox{
(75.2)}\hphantom{0}
\\
&$5\cdot10^3$&
1.30\mbox{ (1.33)}&1.10\mbox{ (1.01)}&3.85\mbox{ (3.77)}&1.21\mbox{ (1.22)}&7.45\mbox{ (6.56)}&2.24\mbox{ (2.20)}&\hphantom{0}16.6\mbox{
(18.1)}\hphantom{0}
\\
&$10^4$&
0.764\mbox{ (0.651)}&0.566\mbox{ (0.474)}&2.50\mbox{ (1.86)}&0.676\mbox{ (0.672)}&4.63\mbox{ (4.03)}&1.01\mbox{ (1.04)}&\hphantom{0}10.1\mbox{ (9.96)}\hphantom{0}\\
\hline
\end{tabular*}
\end{sidewaystable}

In Table \ref{tableunif} we show, for each model with $X$ uniform, the
median (MED) and interquartile range (IQR) of the 200 ISE values
obtained using the LL estimator based on nongrouped data, and, for
several values of $\nu$, the DH and the DM approaches based on data
pooled in groups of size $\nu$; Table~\ref{tablenorm} shows the same
but for $X$ normal. Note that LL cannot be calculated from grouped
data, but we include it to assess the potential loss incurred by
pooling the data. The tables show that for $\nu\leq10$, pooling the
data homogeneously hardly affects the quality of the estimator.
Sometimes, the results are even slightly better with the DH method than
with the LL one. Indeed a careful analysis of the bias and variance of
the various estimators shows that for some curves $p(x)$, grouping
homogeneously can sometimes be slightly beneficial when $\nu$ is small.
(Roughly this is because by grouping a little we lose very little
information, but we increase the number of $Y_j\as$
positive,\vspace*{1pt} which makes the estimation a little easier for
this particular estimator. Theoretical arguments support this
conclusion.) The situation is much less favorable for the DM random
grouping method, whose quality degrades quickly as $\nu$ increases.
Unsurprisingly, DH beat DM systematically, except when $N/\nu$ was
small ($N=1000$ and $\nu =20$), where the $J=50$ grouped observations
did not suffice to estimate very well the curves from models (i) and
(ii).

%s4.2 ###
%s4.2 #&#
\subsection{Real data application}
We also applied our DH method on real data. To make the comparison with
the LL estimator possible, we used data for which we had access to the
entire, nongrouped set of observations $(X_i,Y_i)$. Then we grouped the
data and compared the DH and LL procedures. We used data from the
NHANES study, which are available at
\href{http://www.cdc.gov/nchs/nhanes/nhanes1999-2000/nhanes99\_00.htm}{www.cdc.gov/}
\href{http://www.cdc.gov/nchs/nhanes/nhanes1999-2000/nhanes99\_00.htm}{nchs/nhanes/nhanes1999-2000/nhanes99\_00.htm}. These data
were collected in the US between 1999 and 2000.

As in \citet{DelMei11}, our goal was to estimate two
conditional probabilities: $p_{\mathrm{HBc}}(x)=E(Y_{\mathrm
{HBc}}\mi X=x)$ and \mbox{$p_{\mathrm{CL}}(x)=E(Y_{\mathrm{CL}}\mi X=x)$}, where
$X$ was the age of a patient, $Y_{\mathrm{HBc}}=0$ or $1$ indicating
the absence or presence of antibody to hepatitis B virus core antigen
in the patient's serum or plasma and $Y_{\mathrm{CL}}=0$ or $1$
indicating the absence or presence of genital Chlamydia trachomatis
infection in the urine of the patient. The sample size was $N=7016$
for $\mathrm{HBc}$ and $N=2042$ for $\mathrm{CL}$. The percentage of
$Y_i$'s equal to one was $0.047$ in the HBc case and $0.044$ in the CL
case. See \citet{DelMei11} for more details on these data
and the methods employed to collect them.

%f1 ###
%f1 #&#
\begin{figure}

\includegraphics{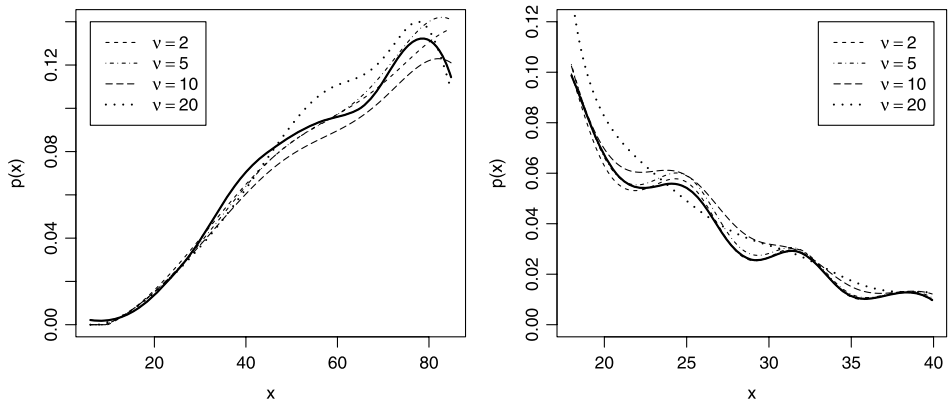}

\caption{NHANES study: DH estimator for $\nu=2$, $5$, $10$ and $20$ and
LL estimator (thick curve) when $Y=Y_{\mathrm{HBc}}$ (left) or
$Y=Y_{\mathrm{CL}}$ (right).}\label{FNhnanes}
\end{figure}

For brevity here we only present the results obtained using our method
by pooling the data homogeneously in groups of equal size $\nu=2$,
$5$, $10$ and $20$. As in the simulations, our DH estimator improved
considerably on the DM method. An illustration of our procedure with a
second covariate is given in Section A.3 of the supplemental article
[\citet{DelHal}].
In Figure~\ref{FNhnanes} we compare DH with LL. All curves were
calculated using our bandwidth procedure described in Section A.2 of
the supplemental article [\citet{DelHal}]. We see that, in
these examples, grouping data in pools of size as large as $\nu=20$
does not dramatically degrade performance.\looseness=1

%s5 ###
%s5 #&#
\section{Generalizations to unequal groups and the multivariate
case}\label{secgener}

Our procedure for estimating $p$ can be extended to the multivariate
setting, where the covariates are random $d$-vectors, and to unequal
group sizes. These extensions can be performed in many different ways,
for example, by binning on each variable, using bins of potentially
different sizes to accommodate different levels of homogeneity.
If we group using bins of equal dimension, then, to a large extent, the
theoretical properties discussed earlier, in the setting of equal-size
groups, continue to hold. To briefly indicate this we give, below,
details of methodology and results in the case of multivariate
histogram binning where, for definiteness, the bin sizes and shapes,
but not the group sizes, are equal. Cases where the bin sizes and
shapes also vary can be treated in a similar manner, provided the
variation is not too great, but since there are so many possibilities
we do not treat those cases here. An approach of this type is discussed
in Section A.3 of the supplemental article [\citet{DelHal}].

In the analysis below we take $\bbX$ to be a $d$-vector, and the
function $p$ to be $d$-variate, where $d\geq1$. We group the data in
bins of equal width, specifically width $(\nu/N)^{1/d}$ along each of
the $d$ coordinate axes, rather than in groups of equal number. In the
theory described below, for notational simplicity, we assume that the
support of the distribution of $\bbX$ contains the cube $\cI
=[0,1]^d$, and we estimate $p$ there. We choose $\nu$ so that
$J=(N/\nu)^{1/d}$ is an integer (on this occasion $\nu$ is not
necessarily an integer itself), and take the bins to be the cubes $\cB
(k_1,\ldots,k_d)$ defined by
\begin{eqnarray*}
\cB(k_1,\ldots,k_d)
&=&\prod_{\ell=1}^d \biggl(\thf(2k_\ell+1) (\nu/N)^{1/d}-\thf
(\nu/N)^{1/d},\\[-6pt]
&&\hspace*{25pt}\thf(2k_\ell+1) (\nu/N)^{1/d}+\thf(\nu/N)^{1/d}\biggr]
,\vspace*{-2pt}
\end{eqnarray*}
where $k_\ell=0,\ldots,J-1$ for $\ell=1,\ldots,d$. In this setting
it is convenient to write the paired data as simply $(\bbX
_1,Y_1),\ldots,(\bbX_N,Y_N)$, where $\bbX_j$ is a $d$-vector and
each $Y_j=0$ or 1, and refer to $\bbX_j$ in terms of the bin in which
it lies, rather than give it a double subscript (as in the notation
$\bbX_{ij}$, where $j$ is the bin index).

Put $b(k_1,\ldots,k_d)=(\thf(2k_1+1) (\nu/N)^{1/d},\ldots,\thf
(2k_d+1) (\nu/N)^{1/d})$, representing the center of the bin $\cB
(k_1,\ldots,k_d)$, define
\[
Z\as(k_1,\ldots,k_d)
=1-\max_{j \dvtx \bbX_j\in\cB(k_1,\ldots,k_d)} Y_j\as\vspace*{-2pt}
\]
and compute $\hmu$ by applying a $d$-variate local polynomial smoother
to the values of $(b(k_1,\ldots,k_d),Z\as(k_1,\ldots,k_d))$,
interpreted as (explanatory variable, response variable) pairs in a
conventional $d$-variate nonparametric regression problem.
To derive an estimator of $p$ from $\hmu$ we take
%
%e5.1 ###
%e5.1 #&#
\begin{equation}\label{52}
\hp({\mathbf x})=1-\hmu({\mathbf x})^{1/m({\mathbf x})},\vspace*{-2pt}
\end{equation}
where $m({\mathbf x})$ denotes the number of data $\bbX_j$ in the
bin containing ${\mathbf x}\in\cI$.\vadjust{\goodbreak}

In developing theoretical properties of this estimator we choose our
regularity conditions to simplify exposition. In particular, we replace
assumptions (S1)--(S4) and (T5) by the following restriction:
{\renewcommand{\thecondition}{U}
\begin{condition}\label{condU}
The nonparametric smoother defined by the estimator at (\ref
{23}) is a standard $d$-variate local linear smoother [see, e.g., \citet{Fan93}], where the kernel $K$, a function of $d$ variables, is a
spherically symmetric, compactly supported, H\"older continuous
probability density, and, for some $\ep>0$, the bandwidth $h$
satisfies $h+(Nh^d)\mo=O(N^{-\ep})$ as $N\rai$.
\end{condition}}

Conditions (T1)--(T4) are replaced by (V1)--(V4) below, and (V5) is
additional:
{\renewcommand{\thecondition}{V}
\begin{condition}\label{condV}

(V1) the distribution of $\bbX$ has a continuous density, $f$,
that is bounded away from zero on an open set $\cJ$ that contains the
cube $\cI=[0,1]^d$;

(V2) the function $p=\delta\pi$ is bounded below 1 uniformly on
$\cI$ and in $N\geq1$;

(V3) the fixed, nonnegative function $\pi$ has two H\"
older-continuous derivatives on $\cJ$;

(V4) for some $\ep>0$, $h+\nu\delta h+(\nu^2/N^{1-\ep}h^d\delta
)\ra0$ as $N\rai$;

(V5) $C_1 (\delta N)^4\leq\nu^{d+4}\leq C_2 N^{d+3}/\delta$
for constants $C_1,C_2>0$.
\end{condition}
%
%th5.1 #&#
\begin{theo}\label{Theorem51} Assume that Conditions \ref{condU} and \ref{condV} hold,
and that $\nu\delta=O(1)$. Then, for each ${\mathbf x}\in\cI$,
%
%e5.2 ###
%e5.2 #&#
\begin{equation}\label{53}
\hp({\mathbf x})=p({\mathbf x})+O_p\{(\delta/Nh^d)\half
+\delta h^2\} .
\end{equation}
\end{theo}

The ``$O_p$'' term on the right-hand side of (\ref{53}) has exactly
the same size as the dominant remainder term, $A({\mathbf x})
V({\mathbf x})+B({\mathbf x})$, on the right-hand side of~(\ref{T7})
in Theorem \ref{Theorem1}, provided of course that we take $d=1$ in
Theorem~\ref{Theorem51}. Refinements given in Theorem \ref{Theorem1}
and in the results in Section \ref{secOP} can also be derived in the
present setting.

Theorem \ref{Theorem51} is proved similarly to Theorem \ref
{Theorem1}, and so is not derived in detail here. The main difference
in the argument comes from incorporating a slightly different
definition of $\hp$, given by (\ref{52}). For example, suppose~$\hp$
is as defined at (\ref{52}), and note that $E(m)=\nu_1+O\{
\nu_1 (\nu_1/N)^2\}$, where $\nu_1({\mathbf x})=\nu
f({\mathbf x})$ and $f$ denotes the density of $\bbX$. Since, in
addition, $m-E(m)=O_p(\nu\half)$, then $m=\nu_1 (1+\De)\mo$ where
$|\De|=O_p\{\nu\mhf+(\nu/N)^2\}$, and, much as in the argument
leading to (\ref{P7}),
%
%e5.3 ###
%e5.3 #&#
\begin{eqnarray}\label{54}
\hp&=&1-\hmu^{1/m}=1-(\hmu^{1/\nu_1})^{1+\De}\nonumber\\
&=&1-[1-p+O_p\{\delta h^2+(\delta/Nh^d)\half
\}]^{1+\De}\nonumber\\[-8pt]\\[-8pt]
&=&1-(1-p) [1+O_p\{\delta h^2+(\delta/Nh^d
)\half+\delta|\De|\}]\nonumber\\
&=&p+O_p[\delta h^2+(\delta/Nh^d)\half
+\delta\{\nu\mhf+(\nu/N)^2\}].\nonumber
\end{eqnarray}
Now, $\delta h^2+(\delta/Nh^d)\half$ is minimized by taking
$h=(N\delta)^{-1/(d+4)}$, and for this choice of $h$ we have
\[
\delta\mo\{\delta h^2+(\delta/Nh^d)\half\}
\asymp(\delta N)^{-2/(d+4)} .
\]
This quantity is not of smaller order than $\nu\mhf+(\nu/N)^2$ if
and only if both $\nu\mhf=O\{(\delta N)^{-\rho}\}$ and $(\nu
/N)^2=O\{(\delta N)^{-\rho}\}$, where $\rho=2/(d+4)$. This is in turn
equivalent to
\[
C_1 (\delta N)^{4/(d+4)}\leq\nu\leq C_2 (N^{d+3}/\delta
)^{1/(d+4)}
\]
for constants $C_1,C_2>0$, which is also equivalent to (V5). Therefore
if (V5) holds, then we can deduce (\ref{53}) from (\ref{54}).

%s6 ###
%s6 #&#
\section{Technical arguments}\label{secP}

%s6.1 ###
%s6.1 #&#
\subsection{\texorpdfstring{Proof of Theorem \protect\ref{Theorem1}}{Proof of Theorem 3.1}}\label{secP1}

Let $D_j$ equal the maximum of $|X_{ij}-\bX_j|$ over $i=1,\ldots,\nu$.
The ratio $\nu/N$ equals the order of magnitude of the expected value
of the width of the group that contains $x\in\cI$, and it can be proved
that
%
%e6.1 #&#
\begin{equation}
\label{P1}
\begin{tabular}{p{300pt}}
for each $\ep>0$, $D_j=O_p(\nu/N^{1-\ep})$ uniformly in $j$ such
that \mbox{$|\bX_j-x|\leq C h$}
and $x\in\cI$.
\end{tabular}\hspace*{-35pt}
\end{equation}
Note that, by (T4), $\nu/N^{1-\ep}\ra0$ for sufficiently small $\ep>0$.

For $k=1,2$ let $p^{(k)}$ be the $k$th derivative of $p$, and put
$p_k=p^{(k)}/\{k! (1-p)\}$. Let $\eta>0$ denote the exponent of
H\"older continuity of $p''$ on $\cI$; see~(T3); that is,
$|p''(x_1)-p''(x_2)|=O(|x_1-x_2|^\eta)$ uniformly in $x_1,x_2\in\cI
$. Then, using~(\ref{P1}) it can be proved that for each $\ep>0$,
%
%e6.2 ###
%e6.2 #&#
\begin{eqnarray}\label{P2}
E(Z_j\as\mi\cX)&=&\prod_{i=1}^\nu\{1-p(X_{ij})\}\nonumber\\
&=&\{1-p(\bX_j)\}^\nu\prod_{i=1}^\nu\{1-p_1(\bX_j) (X_{ij}
-\bX_j)
+O_p(\delta D_j^2)\}\nonumber\\
&=&\{1-p(\bX_j)\}^\nu\prod_{i=1}^\nu\exp\{-p_1(\bX_j)
(X_{ij}-\bX_j)
+O_p(\delta D_j^2)\}\nonumber\\[-8pt]\\[-8pt]
&=&\{1-p(\bX_j)\}^\nu\exp\Biggl\{-\sum_{i=1}^\nu p_1(\bX_j)
(X_{ij}-\bX_j)
+O_p(\nu\delta D_j^2)\Biggr\}\nonumber\\
&=&\{1-p(\bX_j)\}^\nu\exp\{O_p(\nu\delta D_j^2
)\}\nonumber\\
&=&\{1-p(\bX_j)\}^\nu\{1+O_p(\nu^3 \delta/N^{2-\ep
})\} ,\nonumber
\end{eqnarray}
uniformly in the sense of (\ref{P1}) and for each $\ep>0$.
[Assumption (T4) implies that $\nu^3\delta/N^{2-\ep}\ra0$ for some\vadjust{\goodbreak}
$\ep>0$.] Observe too that, uniformly in the same sense,
%
%e6.3 ###
%e6.3 #&#
\begin{eqnarray}\label{P3}
&&\{1-p(\bX_j)\}^\nu\nonumber\\
&&\qquad=\{1-p(x)\}^\nu\{1-p_1(x) (\bX_j-x)\nonumber\\
&&\qquad\quad\hspace*{55.5pt}{}-p_2(x) (\bX
_j-x)^2
+O_p(\delta h^{2+\eta})\}^\nu\nonumber\\[-8pt]\\[-8pt]
&&\qquad=\{1-p(x)\}^\nu\bigl[1-\nu p_1(x) (\bX_j-x)\nonumber\\
&&\qquad\quad\hspace*{55.5pt}{}+\bigl\{\tfrac12\nu(\nu-1) p_1(x)^2
-\nu p_2(x)\bigr\} (\bX_j-x)^2\nonumber\\
&&\qquad\quad\hspace*{131pt}{}+O_p(\nu\delta h^{2+\eta}
+\nu^3 \delta^3 h^3)\bigr] ,\nonumber
\end{eqnarray}
again uniformly in the sense of (\ref{P1}). [Note that, by (T4),
$\nu\delta h\ra0$.] Combining (\ref{23}), (T4), (S1), (S2), (S4),
(\ref{P2}) and (\ref{P3}) we deduce that, for each $\ep>0$ and
each $x\in\cI$,
\begin{eqnarray*}
\tmu(x)
&\equiv& E\{\hmu(x)\mi\cX\}
={\sumj w_j(x) E(Z_j\as\mi\cX)\Big/\sumj w_j(x)}\\
&=&\{1-p(x)\}^\nu\biggl\{1+\biggl[\thf\nu(\nu-1) p_1(x)^2
-\nu p_2(x)\biggr] {\sumj w_j(x) (\bX_j-x)^2\over\sumj
w_j(x)}\\
&&\hspace*{135.4pt}{}
+O_p(\nu\delta h^{2+\eta}+\nu^3 \delta^3 h^3+\nu^3
\delta N^{\ep-2})\biggr\}\\
&=&\{1-p(x)\}^\nu\biggl[1+h^2 \biggl\{\thf\nu(\nu-1) p_1(x)^2
-\nu p_2(x)\biggr\} b(x)\\
&&\hspace*{134pt}{}
+o_p(\nu\delta h^2+\nu^3 \delta N^{\ep-2})\biggr] ,
\end{eqnarray*}
whence, for all $\ep>0$,
%
%e6.4 ###
%e6.4 #&#
\begin{eqnarray}\label{P4}
\tmu(x)^{1/\nu}&=&\{1-p(x)\}\bigl[1-h^2\bigl\{p_2(x)-\tfrac12(\nu
-1)p_1(x)^2\bigr\} b(x)\nonumber\\[-8pt]\\[-8pt]
&&\qquad\quad\hspace*{85pt}{}+o_p(\delta h^2+\nu^2\delta N^{\ep
-2})\bigr] ,\nonumber
\end{eqnarray}
uniformly in $x\in\cI$. Hence, defining
%
%e6.5 ###
%e6.5 #&#
\begin{equation}\label{P5}
\De(x)=\hmu(x)-\tmu(x)
={\sumj w_j(x) \{Z_j\as-E(Z_j\as\mi\cX)\}\Big/\sumj w_j(x)},
\end{equation}
noting that $1-p$ is bounded away from zero [see (T2)], and taking the
argument of the functions below to equal the specific point $x$
referred to in~(\ref{T7}), we deduce that
%
%e6.7 ###
%e6.6 ###
%e6.6 #&#
%e6.7 #&#
\begin{eqnarray}
\hp&=&1-\hmu^{1/\nu}=1-(\tmu+\De)^{1/\nu}\nonumber\\[-1pt]
&=&1-\bigl(\tmu^{1/\nu}+\nu\mo\tmu^{-(\nu-1)/\nu} \De\bigr)
+O_p\bigl(\nu\mo\tmu^{-(2\nu-1)/\nu} \De^2\bigr)\nonumber\\[-1pt]
&=&1-(1-p) \bigl[1-h^2 \bigl\{p_2-\tfrac12(\nu-1) p_1^2\bigr\} b
+o_p(\delta h^2+\nu^2 \delta N^{\ep-2})\bigr]\nonumber
\\[-1pt]
&&{}
-\nu\mo\tmu^{-(\nu-1)/\nu} \De
+O_p\bigl(\nu\mo\tmu^{-(2\nu-1)/\nu} \De^2\bigr)\nonumber\\[-1pt]
&=&p+(1-p) \bigl[h^2 \bigl\{p_2-\tfrac12(\nu-1) p_1^2\bigr\} b
+o_p(\delta h^2+\nu^2 \delta N^{\ep-2})\bigr]\nonumber
\\[-1pt]
\label{P6}
&&{}
-\{1+o_p(1)\} \nu\mo(1-p)^{-(\nu-1)} \De
+O_p\bigl\{\nu\mo(1-p)^{-(2\nu-1)}\De^2\bigr\} \\[-1pt]
&=&p+(1-p) \bigl[h^2 \bigl\{p_2-\tfrac12(\nu-1) p_1^2\bigr\} b
+o_p(\delta h^2)\bigr]\nonumber\\[-1pt]
\label{P7}
&&{}
-\{1+o_p(1)\} \nu\mo(1-p)^{-(\nu-1)} \De,
\end{eqnarray}
where (\ref{P6}) holds without the assumption $\nu\delta=O(1)$
[it holds under either that condition or (\ref{T22})], but (\ref
{P7}) requires $\nu\delta=O(1)$. Note that, by (T4), $\nu
/N^{1-\ep
}h\ra0$ for some $\ep>0$, and so $\nu^2\delta N^{2\ep-2}/(\delta
h^2)=(\nu/N^{1-\ep}h)^2\ra0$. Additionally, it will follow from
(\ref{P9}) below that, when $\delta=O(1)$, $\De=O_p\{(\nu
^2\delta/Nh)\half\}$, and by~(T4), $\delta/Nh\ra0$, so $\De
=o_p(1)$. The identity leading from (\ref{P6}) to (\ref{P7})
follows from this property.

Observe that, by (\ref{P2}) and (\ref{P3}), $E(Z_j\as\mi\cX
)=\{1+o_p(1)\} \{1-p(x)\}^\nu$ and
\[
1-E(Z_j\as\mi\cX)
=1-\{1-p(x)\}^\nu
+O_p[\{1-p(x)\}^\nu(\nu\delta h+\nu^3 \delta
N^{\ep-2})] ,
\]
uniformly in $j$ such that $|\bX_j-x|\leq C h$, where $C$ is as in
(T5), and moreover,
\[
\var(\De\mi\cX)
={\sumj w_j^2 \var(Z_j\as\mi\cX)\over(\sumj w_j)^2}
={\sumj w_j^2 E(Z_j\as\mi\cX) \{1-E(Z_j\as\mi\cX)\}\over(\sumj
w_j)^2} .
\]
[Here and in (\ref{P8})--(\ref{P10}) the argument of the
functions is the point $x$ in~(\ref{T7}).] Therefore, by (S3),
%
%e6.8 ###
%e6.8 #&#
\begin{eqnarray}\label{P8}
\var(\De\mi\cX)
&=&\{1+o_p(1)\} (\nu/Nh) (1-p)^\nu\{1-(1-p)^\nu\}
v\nonumber\\[-8pt]\\[-8pt]
&&{}+O_p[(\nu/Nh) (\nu\delta h+\nu^3 \delta N^{\ep
-2})] .\nonumber
\end{eqnarray}
Properties (T4) and (\ref{P8}), and Lyapounov's central limit
theorem (see the next paragraph for details), imply that when $\nu
\delta=O(1)$ and $\pi(x)>0$ [the latter is assumed here and below;
the proof when $\pi(x)=0$ is simpler], we can write
%
%e6.9 ###
%e6.9 #&#
\begin{eqnarray}\label{P9}
\De&=&\bigl((\nu/Nh)(1-p)^\nu\{1-(1-p)^\nu\} v\nonumber\\[-1pt]
&&\hspace*{2pt}{}+O_p[(\nu/Nh)(\nu\delta h+\nu^3\delta N^{\ep-2})
]\bigr)\half V_4\\[-1pt]
&=&\{1+o_p(1)\}
[(\nu/Nh) (1-p)^\nu\{1-(1-p)^\nu\} v]\half
V_4 ,\nonumber
\end{eqnarray}
where the second identity follows from the fact that $h+\nu^2 N^{\ep
-2}\ra0$ for some $\ep>0$ [see (T4)], and $V_4$ denotes a random
variable that is asymptotically distributed as normal N$(0,1)$. This
result and (\ref{P7}) imply that
%
%e6.10 ###
%e6.10 #&#
\begin{eqnarray}\label{P10}
\hp&=&p+(1-p) \bigl[h^2 \bigl\{p_2-\tfrac12(\nu-1) p_1^2\bigr\} b
+o_p(\delta h^2)\bigr]\nonumber\\[-9pt]\\[-9pt]
&&{}
-\{1+o_p(1)\} [(\nu Nh)\mo(1-p)^{2-\nu}
\{1-(1-p)^\nu\} v]\half V_4 .\nonumber
\end{eqnarray}
Result (\ref{T7}) follows from (\ref{P10}).\vadjust{\goodbreak}

When applying a generalized from of Lyapounov's theorem to establish
a~central limit theorem for $\De$, conditional on $\cX$, we should, in
view of (S4), prove that for some integer $k>2$,
$
[(\nu/ Nh) (1-p)^\nu\{1-(1-p)^\nu\} v]^{-k/2}
(\nu/\break Nh)^{k-1}\ra0 .
$
When $\nu\delta=O(1)$ this is equivalent to $(\delta/Nh)^{-k/2}(\nu
/Nh)^{k-1}\ra0$, and hence to $(Nh/\nu)^2 (\nu^2/Nh\delta)^k\ra
0$; call this result (R). Now, (T4) ensures that for some $\ep>0$,
$\nu^2/N^{1-\ep}h\delta\ra0$. Therefore (R) holds for all
sufficiently large $k$.

Next we outline the derivation of (\ref{T10}). It can be proved
from (\ref{T6}) that if $C_1>0$ is given, if $C_2=C_2(C_1)>0$ is
chosen sufficiently large, if $\cI_N$ is a~regular grid of $n^{C_2}$
points in $\cI$ and if, for each $x\in\cI$, we define $x_N$ to be
the point in $\cI_N$ nearest to $x$, then
%
%e6.11 ###
%e6.11 #&#
\begin{equation}\label{P11}
P\Bigl\{{\sup_{x\in\cI}} |\De(x)-\De(x_N)|\leq N^{-C_1}\Bigr\}\ra1 .
\end{equation}
Note that, by (T4),
applying (S3), (S4), Rosenthal's and Markov's inequalities, we can
prove that, for each $C,\ep>0$,
$
\sup_{x\in\cI} P\{|\De(x)|>N^\ep(\nu^2\delta/Nh)\half
\mi\cX\}
=O_p(N^{-C}) .
$
It follows that, for all $C,\ep>0$,
%
%e6.12 ###
%e6.12 #&#
\begin{equation}\label{P13}
P\Bigl\{{\sup_{x\in\cI_N}} |\De(x)|>N^\ep(\nu^2\delta/Nh)\half
\Bigmi\cX\Bigr\}
=O(N^{-C}) .
\end{equation}
Together (\ref{P11}) and (\ref{P13}) imply that, for each
$C,\ep>0$,
%
%e6.13 ###
%e6.13 #&#
\begin{equation}\label{P14}
P\Bigl\{{\sup_{x\in\cI}} |\De(x)|>N^\ep(\nu^2\delta/Nh)\half
\Bigr\}
\ra0 .
\end{equation}

Results (\ref{P4}) (which holds uniformly in $x\in\cI$) and
(\ref{P14}) imply that (\ref{P7}) holds uniformly in $x\in\cI
$. Hence,
%
%e6.14 ###
%e6.14 #&#
\begin{eqnarray}\label{P15}
\inti(\hp-p)^2
&=&\inti B^2+\inti\bigl\{\nu\mo(1-p)^{-(\nu-1)} \De\bigr\}
^2\nonumber\\
&&{}-2\inti B \bigl\{\nu\mo(1-p)^{-(\nu-1)} \De\bigr\}\\
&&{}+o_p\biggl\{(\delta h^2)^2
+\inti(\De/\nu)^2\biggr\}.\nonumber
\end{eqnarray}
Conditional on $\cX$ the random variable $\De$, at (\ref{P5}),
equals a sum of independent random variables with zero means, and using
that property, Condition~\ref{condS} (which, for this part of the theorem, holds
uniformly in $x\in\cI$) and (\ref{T7}), it can be proved that
%
%e6.17 ###
%e6.16 ###
%e6.15 ###
%e6.15 #&#
%e6.16 #&#
%e6.17 #&#
\begin{eqnarray}
\label{P16}
E\biggl[\inti\bigl\{\nu\mo(1-p)^{-(\nu-1)} \De\bigr\}^2\Bigm|
\cX\biggr]
&=&\inti A^2+o_p(\delta/Nh) ,\\
\label{P17}
\var\biggl[\inti\bigl\{\nu\mo(1-p)^{-(\nu-1)} \De\bigr\}
^2\Bigm|\cX\biggr]
&=&o_p\{(\delta/Nh)^2\} ,\\
\label{P18}
\var\biggl[\inti B \bigl\{\nu\mo(1-p)^{-(\nu-1)} \De\bigr\}
\Bigm|\cX\biggr]
&=&o_p\{(\delta/Nh)^2+(\delta h^2)^4\} .
\end{eqnarray}
Result (\ref{P16}) follows from (\ref{P8}). To derive (\ref
{P17}), note that by (\ref{P4}) we have, uniformly in
$x_1,x_2\in
\cI$,
%
%e6.18 ###
%e6.18 #&#
\begin{eqnarray}\label{A}
E\{\De(x_1)^2 \De(x_2)^2\mi\cX\}
&=&E\{\De(x_1)^2\mi\cX\} E\{\De(x_2)^2\mi
\cX\}\nonumber\\[-8pt]\\[-8pt]
&&{}+O_p\{t_1(x_1,x_2)\} ,\nonumber
\end{eqnarray}
where
\begin{eqnarray*}
t_1(x_1,x_2)
&=&{\sumj w_j(x_1)^2 w_j(x_2)^2 E[\{Z_j\as-E(Z_j\as\mi\cX)\}
^4\mi\cX]
\over\{\sumj w_j(x_1)\}^2 \{\sumj w_j(x_2)\}^2}
=O_p\{t_2(x_1,x_2)\} ,\\
t_2(x_1,x_2)&=&{\sumj w_j(x_1)^2 w_j(x_2)^2
\var(Z_j\as\mi\cX)
\over\{\sumj w_j(x_1)\}^2 \{\sumj w_j(x_2)\}^2}
=O_p\biggl[{\nu\delta\sumj w_j(x_1)^4\over\{\sumj w_j(x_1)\}
^4}\biggr]\\
&=&O_p\biggl\{\biggl({\nu\delta\over Nh}\biggr)^{ 2} \biggl({\nu^2\over
Nh\delta}\biggr)\biggr\}
=o_p\biggl\{\biggl({\nu\delta\over Nh}\biggr)^{ 2}\biggr\} ,
\end{eqnarray*}
again uniformly in $x_1,x_2\in\cI$. [The last and second-last
identities here follow from (T4) and (S4), resp.] Noting these
bounds, defining $\xi_1\equiv\{\nu\mo(1-p)^{-(\nu-1)}\}^2$ and
integrating (\ref{A}) over $x_1,x_2\in\cI$, we deduce that
\[
E\biggl\{\int\xi_1(x) \De(x)^2 \,dx\Bigm|\cX\biggr\}^2
=\biggl[\int\xi_1(x) E\{\De(x)^2\mi\cX\} \,dx\biggr]^2
+o_p\{(\nu\delta/Nh)^2\} ,
\]
which implies (\ref{P17}).

%We alo need

To derive (\ref{P18}), define $\xi_2=B \xi_1$ and $e_j=E[\{
Z_j\as-E(Z_j\as\mi\cX)\}^2\mi\cX]$, write~$M$ for the left-hand
side of (\ref{P18}), and note that
\[
M=\inti\inti\xi_2(x_1) \xi_2(x_2)
{\sumj w_j(x_1) w_j(x_2) e_j
\over\{\sumj w_j(x_1)\} \{\sumj w_j(x_2)\}} \,dx_1 \,dx_2 .
\]
In view of (T5), $w_j(x)=0$ if $|\bX_j-x|>C h$, and so the series in
the numerator inside the integrand can be confined to indices $j$ for
which both $|\bX_j-x_1|\leq C h$ and $|\bX_j-x_2|\leq C h$.
Therefore the integrand equals zero unless $|x_1-x_2|\leq2 C h$.
Hence, defining $J(x_1,x_2)=1$ if $|x_1-x_2|\leq2 C h$, and
$J(x_1,x_2)=0$ otherwise, using the Cauchy--Schwarz inequality to
derive both the inequalities below and writing $\|\cI\|$ for the
length of the interval $\cI$, we have
%
%e6.19 ###
%e6.19 #&#
\begin{eqnarray}\label{B}
M&\leq&\inti\inti J(x_1,x_2)\xi_2(x_1)\xi_2(x_2)\nonumber\\
&&\hspace*{23pt}{}\times\Biggl[\prod_{k=1}^2{\sumj w_j(x_k)^2 e_j
\Big/\biggl\{\sumj w_j(x_k)\biggr\}^2}\Biggr]\half \,dx_1\, dx_2
\nonumber\\[-8pt]\\[-8pt]
&=& \inti\inti J(x_1,x_2) \xi_2(x_1) \xi_2(x_2)
\Biggl[\prod_{k=1}^2 \var\{\De(x_k)\mi\cX\}\Biggr]\half \,dx_1\,
dx_2\nonumber\\[-2pt]
&\leq&\|\cI\| \Biggl(\inti\inti J(x_1,x_2)
\Biggl[\prod_{k=1}^2 \xi_2(x_k)^2 \var\{\De(x_k)\mi\cX\}
\Biggr]\,
dx_1 \,dx_2\Biggr)^{ 1/2}.\nonumber
\end{eqnarray}
Using (\ref{P8}) show that $\var\{\De(x_k)\mi\cX\}=O(\nu
^2\delta/Nh)$, uniformly in $x_k\in\cI$, noting that $B=O(\delta
h^2)$ uniformly in $x\in\cI$ [the bound at (\ref{T10}) holds uniformly in
the argument of $B$] and observing that $\xi_1(x)=O(\nu\mt)$
uniformly in $x\in\cI$, whence it follows from the bound for $B$ that
$\xi_2(x)=O(\delta h^2 \nu\mt)$ uniformly in $x\in\cI$, we
deduce from (\ref{B}) that
%
%e6.20 ###
%e6.20 #&#
\begin{eqnarray}\label{C}
M&=&O_p[(\nu^2\delta/Nh) \{(\delta
h^2 \nu\mt)\}^2]
\biggl(\inti\inti J(x_1,x_2) \,dx_1 \,dx_2\biggr)^{
1/2}\nonumber\\[-8pt]\\[-8pt]
&=&O_p(\delta^2 h^{7/2}/N)=o_p\{(\delta
/Nh)^2+(\delta h^2)^4\} .\nonumber
\end{eqnarray}
Result (\ref{P18}) follows directly from (\ref{C}).

%For closely related results in regression problems, see for example
%Hall (1984).

%estimators of regression functions. {\sl Ann. Statist.} {\bf12},
%241--260.

%s6.2 ###
%s6.2 #&#
\subsection{\texorpdfstring{Proof of Theorem \protect\ref{Theorem2}}{Proof of Theorem 3.2}}\label{secP2}

The proof is similar to that of the first part of Theorem
\ref{Theorem1}, the main difference occurring at the point at which the
remainder term, $O_p(R)$ where $R=\nu\mo(1-p)^{-(2\nu-1)} \De
^2$, in
(\ref{P6}), is shown to be negligible relative to the term
$\nu\mo(1-p)^{-(\nu-1)} \De$ there. It suffices to prove that
$(1-p)^{-\nu} \De\ra0$ in probability, or equivalently, in view of
(\ref{P9}), that $(\nu/Nh) (1-p)^{-\nu}\ra0$. However, the latter
result is ensured by (\ref{T22}).

%s6.3 ###
%s6.3 #&#
\subsection{\texorpdfstring{Proof of Theorem \protect\ref{Theorem3}}{Proof of Theorem 3.3}}\label{secP3}
Without loss of generality, the point $x$ in~(\ref{T27}) is $x=0$.
Recall\vspace*{1pt} that $p^0\equiv\delta$, and take $p^1(u)=\delta\{1+h^2
\psi(u/h)\}$, where~$\psi$ is bounded and has two bounded derivatives
on the real line, is supported on $[-\thf,\thf]$ and satisfies $\psi
(0)\neq0$. The respective functions $\pi^0\equiv1$ and $\pi
^1(u)=1+h^2 \psi(u/h)$ satisfy~(\ref{T25}). [The quantity
$h=h(N)>0$ here is not a bandwidth, but converges to 0 as $N\rai$.]
Therefore, $p^0(u)=p^1(u)$ except when $u\in(-\thf h,\thf h)$. We
assume that $\nu\delta\rai$ as $N\rai$, and consider the problem of
discriminating between $p^0$ and $p^1$ using the data pairs $(\cX
_j,Y_j\as)$.

Without loss of generality, we confine attention to those pairs $(\cX
_j,Y_j\as)$ for which $\cX_j$ is wholly contained in $[-\thf h,\thf
h]$. Pairs for which $\cX_j$ has no intersection with $[-\thf
h,\thf h]$ convey no information for discriminating between $p^0$
and $p^1$, and it is readily proved that including pairs for which $\cX
_j$ overlaps the boundary does not affect the results we derive below.
In a slight abuse of notation we shall take the integers $j$ for which
$\cX_j\subseteq[-\thf h,\thf h]$ to be $1,\ldots,m$, where
$m=hN/\nu+o_P(1)$ and is assumed to be an integer.

The likelihood\vspace*{-2pt} of the data pairs $(\cX_j,Y_j\as)$ for
$1\leq j\leq m$, conditional on $\cX=\{X_1,\ldots,X_N\}$, is $
\prod_{j=1}^m P_j^{Y_j\as} (1-P_j)^{1-Y_j\as} $ where $
P_j=P(Y_j\as=1\mi\cX)=1-\prod_{i=1}^\nu\{1-p(X_{ij})\} . $ Let
$P_j^0$\vadjust{\goodbreak}
and $P_j^1$ denote the versions of $P_j$ when $p=p^0$ and $p=p^1$,
respectively. Also, let $\Theta_j^+=P_j^1/P_j^0$ and $\Theta
_j^-=(1-P_j^1)/(1-P_j^0)$. In this notation the log-likelihood ratio
statistic is given by
%
%e6.21 ###
%e6.21 #&#
\begin{eqnarray}\label{P19}
L&=&\sum_{j=1}^m \{Y_j\as\log(\Theta_j^+)
+(1-Y_j\as) \log(\Theta_j^-)\}\nonumber\\[-8pt]\\[-8pt]
&=&\sum_{j=1}^m (1-Y_j\as) \log(\Theta_j^-/\Theta_j^+)
+\sum_{j=1}^m \log(\Theta_j^+) \nonumber
\end{eqnarray}
and therefore,
$E(L\mi\cX)
=\sum_{j=1}^m (1-P_j) \log(\Theta_j^-/\Theta_j^+)
+\sum_{j=1}^m \log(\Theta_j^+) $,
$\var(L\mi\cX)
=\sum_{j=1}^m P_j (1-P_j)
\{\log(\Theta_j^-/\Theta_j^+)\}^2 .
$
Writing $E^0$ and $\var^0$ to denote expectation and variance when
$p=p^0$, we deduce that
%
%e6.23 ###
%e6.22 ###
%e6.22 #&#
%e6.23 #&#
\begin{eqnarray}
\label{P20}
E^0(L\mi\cX)
&=&(1-\delta)^\nu\sum_{j=1}^m \log(\Theta_j^-/\Theta_j^+)
+\sum_{j=1}^m \log(\Theta_j^+) ,\\
\label{P21}
\var^0(L\mi\cX)
&=&(1-\delta)^\nu\{1-(1-\delta)^\nu\} \sum_{j=1}^m
\{\log(\Theta_j^-/\Theta_j^+)\}^2 .
\end{eqnarray}

Assume for the time being that
%
%e6.24 ###
%e6.24 #&#
\begin{equation}\label{P22}
\nu\delta h^2\ra0
\end{equation}
as $N\rai$, and observe that, since $1-P_j^0=(1-\delta)^\nu$, then
%
%e6.25 ###
%e6.25 #&#
\begin{eqnarray}\label{P23}
\Theta_j^-
&=&(1-P_j^0)\mo\prod_{i=1}^\nu[1-\delta\{
1+h^2 \psi(X_{ij}/h)\}]\nonumber\\[-8pt]\\[-8pt]
&=&\prod_{i=1}^\nu\biggl\{1-{\delta\over1-\delta} h^2 \psi
(X_{ij}/h)\biggr\}
=1-\rho h^2 S_j+R_j ,\nonumber
\end{eqnarray}
where $\rho=\delta/(1-\delta)$, $S_j=\sumi\psi(X_{ij}/h)$ and
$R_j=O_p(\nu\rho^2 h^4)$ uniformly in $1\leq j\leq m$. [We used
(\ref{P22}) to derive the last identity in (\ref{P23}). To
obtain uniformity in the bound for $R_j$, and in later bounds, we used
the fact that~$\psi$ is bounded.] Hence,
\[
\log(\Theta_j^-)
=-\bigl\{\rho h^2 S_j-R_j+\tfrac12(\rho h^2 S_j-R_j)^2
+ \tfrac13 (\rho h^2 S_j-R_j)^3-\cdots\bigr\}.
\]
Similarly, since
\begin{eqnarray*}
P_j^1&=&1-(1-P_j^1)\\
&=&1-(1-P_j^0) \prod_{i=1}^\nu
\biggl\{1-{\delta\over1-\delta} h^2 \psi(X_{ij}/h)\biggr\}\\
&=&1-(1-P_j^0) (1-\rho h^2 S_j+R_j)\\
&=&P_j^0+(1-P_j^0) (\rho h^2 S_j-R_j) ,
\end{eqnarray*}
then
\begin{eqnarray*}
\log(\Theta_j^+)
&=&\log\{1+(\rho h^2 S_j-R_j)(1-P_j^0)/P_j^0\}
\\[-1pt]
&=&{(1-\delta)^\nu\over1-(1-\delta)^\nu} (\rho h^2
S_j-R_j)
-\thf(1-\delta)^{2\nu} (\rho h^2 S_j-R_j)^2
\\[-1pt]
&&{}
+O_p\{(1-\delta)^{3\nu} (\nu\rho h^2)^2\} ,
\end{eqnarray*}
uniformly in $1\leq j\leq m$. It follows that
%
%e6.26 ###
%e6.26 #&#
\begin{eqnarray}\label{P24}
&&\log(\Theta_j^-/\Theta_j^+)\nonumber\\[-1pt]
&&\qquad=-\biggl\{(\rho h^2 S_j-R_j)
+\frac12(\rho h^2 S_j-R_j)^2
+{ \frac13} (\rho h^2 S_j-R_j)^3
+\cdots\biggr\}\nonumber\\[-1pt]
&&\qquad\quad{}
-{(1-\delta)^\nu\over1-(1-\delta)^\nu} (\rho h^2
S_j-R_j)
+\frac12(1-\delta)^{2\nu} (\rho h^2 S_j
-R_j)^2\\[-1pt]
&&\qquad\quad{}+O_p\{(1-\delta)^{3\nu} (\nu\rho h^2
)^2\}\nonumber\\[-1pt]
&&\qquad=-\rho h^2 S_j+O_p\{\nu\rho^2 h^4
+(1-\delta)^\nu\nu\rho h^2\} ,\nonumber\\[-1pt]
%
%e6.27 ###
%e6.27 #&#
\label{P25}
&&(1-\delta)^\nu\log(\Theta_j^-/\Theta_j^+)
+\log(\Theta_j^+)\nonumber\\[-1pt]
&&\qquad=-(1-\delta)^\nu\biggl\{(\rho h^2 S_j-R_j)
+\frac12(\rho h^2 S_j-R_j)^2
+{ \frac13} (\rho h^2 S_j-R_j)^3
+\cdots\biggr\}\hspace*{-4pt}\nonumber\\[-1pt]
&&\qquad\quad{}
-{(1-\delta)^{2\nu}\over1-(1-\delta)^\nu} (\rho h^2
S_j-R_j)\nonumber\\[-1pt]
&&\qquad\quad{}
+{(1-\delta)^\nu\over1-(1-\delta)^\nu} (\rho h^2
S_j-R_j)
-\frac12(1-\delta)^{2\nu} (\rho h^2 S_j-R_j)^2
\nonumber\\[-9pt]\\[-9pt]
&&\qquad\quad{}
+O_p\{(1-\delta)^{3\nu} (\nu\rho h^2)^2\}
\nonumber\\[-1pt]
&&\qquad=-(1-\delta)^\nu\biggl\{\frac12(\rho h^2 S_j-R_j)^2
+{ \frac13} (\rho h^2 S_j-R_j)^3+\cdots
\biggr\}\nonumber\\[-1pt]
&&\qquad\quad{}
-\frac12(1-\delta)^{2\nu} (\rho h^2 S_j-R_j)^2
+O_p\{(1-\delta)^{3\nu} (\nu\rho h^2)^2\}
\nonumber\\[-1pt]
&&\qquad=-\frac12(1-\delta)^\nu(\rho h^2 S_j)^2
+o_p\{(1-\delta)^\nu(\nu\rho h^2)^2\} ,\nonumber
\end{eqnarray}
uniformly in $1\leq j\leq m$. Using (\ref{P20}), (\ref{P21}),
(\ref{P24}) and (\ref{P25}) we deduce that
\begin{eqnarray*}
E^0(L\mi\cX)
&=&-\thf(1-\delta)^\nu(\rho h^2)^2 \sum_{j=1}^m S_j^2
+o_p\{m (1-\delta)^\nu(\nu\rho h^2)^2\}
,\\[-1pt]
\var^0(L\mi\cX)
&=&\{1+o_p(1)\} (1-\delta)^\nu(\rho h^2)^2 \sum
_{j=1}^m S_j^2
+o_p\{m (1-\delta)^\nu(\nu\rho h^2)^2\} .
\end{eqnarray*}

Choose $h$ so that
%
%e6.28 ###
%e6.28 #&#
\begin{equation}\label{P26}
\mbox{the squared mean and the variance are of the same order}.
\end{equation}
In particular, take
$
\{m (1-\delta)^\nu(\nu\rho h^2)^2\}^2
=C_1 m (1-\delta)^\nu(\nu\rho h^2)^2 ,
$
and hence
%
%e6.29 ###
%e6.29 #&#
\begin{equation}\label{P27}
m (1-\delta)^\nu(\nu\rho h^2)^2=C_2+o_P(1)
\end{equation}
or equivalently, using the fact that $m=Nh/\nu+o_P(1)$,
%
%e6.30 ###
%e6.30 #&#
\begin{equation}\label{P28}
h=C_3 \{(N \nu\rho^2)\mo(1-\delta)^{-\nu}\}^{1/5},
\end{equation}
where $C_1,C_2,C_3$ are positive constants; $C_3$ can be chosen
arbitrarily. It follows that
%
%e6.31 ###
%e6.31 #&#
\begin{equation}\label{P29}
\rho h^2
=C_3^2 (\rho/N^2)^{1/5} \nu^{-2/5} \la_N^2,
\end{equation}
where $\la_N^5=(1-\delta)^{-\nu}$. If $h$ is given by (\ref
{P29}), then $\nu\rho h^2=C_3^2 (\nu^3\rho/N^2)^{1/5} \la_N^2$
and therefore (\ref{P22}) follows from (\ref{T26}).

It can be shown that, conditional on the explanatory variables, the
log-likeli\-hood ratio $L$, centred at the conditional mean and
variance, is asymptotically normally distributed with zero mean and
unit variance. (We shall give a proof below.) Therefore by taking
$C_1$, and hence $C_3$, sufficiently small, we can ensure that: (i)~The
probability of discriminating between $p^0$ and $p^1$, when $p=p^0$,
is bounded below~1 as $N\rai$. [This follows from (\ref{P26}).]
Similarly it can be proved that: (ii) The probability of discriminating
between~$p^0$ and~$p^1$, when $p=p^1$, is bounded below 1. Consider
the assertion: (iii)~${\check p}(0)-p(0)$ converges in probability
to 0, along a subsequence, at a~strictly faster rate than $h^2$.
If (iii) is true, then the error rate of the classifier which asserts
that $p=p^0$ if ${\check p}(0)$ is closer to $p(0)$ than to $p^1(0)$,
and $p=p^1$ otherwise, and converges to 0 as $N\rai$. However,
properties (i) and (ii) show that even the optimal classifier, based on
the likelihood ratio rule, does not enjoy this degree of accuracy, and
so (iii) must be false. This proves~(\ref{T27}).

Finally we derive the asymptotic normality of $L$ claimed in the
previous paragraph. We do this using Lindeberg's central limit theorem,
as follows. In view of the definition of $L$ at (\ref{P19}) it is
enough to prove that for each $\eta>0$,\looseness=-1
%
%e6.32 ###
%e6.32 #&#
\begin{eqnarray}\label{P30}\quad
S_{N1}(\eta)&\equiv&\si(\cX)\mt\sum_{j=1}^m
E^0\bigl[\bigl|Y_j\as-E(Y_j\as\mi\cX)\bigr|^2 \si_j(\cX
)^2\nonumber\\[-1pt]
&&\hspace*{70pt}{}
\times I\bigl\{\bigl|Y_j\as-E(Y_j\as\mi\cX)\bigr| \si_j(\cX)
>\eta\si(\cX)\bigr\}\Bigmi\cX\bigr]\\[-1pt]
&\ra&0\nonumber
\end{eqnarray}\looseness=0
in probability, where we define
%
%e6.34 ###
%e6.33 ###
%e6.33 #&#
%e6.34 #&#
\begin{eqnarray}\hspace*{28pt}
\label{P31}
\si_j(\cX)^2
&=&\{\log(\Theta_j^-/\Theta_j^+)\}^2
=(\rho h^2S_j)^2+o_p\{(\nu\delta h^2)^2\} ,\\[-1pt]
\label{P32}
\si(\cX)^2&=&\sum_{j=1}^m \var^0(Y_j\as\mi\cX)
\si_j(\cX)^2
=\{1+o_p(1)\} C_4 m (\nu\rho h^2)^2 (1-\delta)^\nu\vadjust{\goodbreak}
\end{eqnarray}
with $C_4>0$ and (\ref{P31}) holding uniformly in $j$. [We used
(\ref{P24}) to obtain the second identities in each of (\ref
{P31}) and (\ref{P32}).] Since $m=hN/\nu+o_P(1)$, then, by
(\ref{P28}) and (\ref{P32}), $\si(\cX)^2\ra C_5$ in
probability, where $C_5>0$. Hence, by (\ref{P30}), with probability
converging to 1 as $N\rai$,
\begin{eqnarray*}
C_6 S_{N1}(\eta)&\leq& S_{N2}(\eta)\\[-3pt]
&\equiv&\sum_{j=1}^m
E^0\bigl[\bigl|Y_j\as-E(Y_j\as\mi\cX)\bigr|^2 \si_j(\cX)^2\\[-3pt]
&&\hspace*{33pt}{}
\times I\bigl\{\bigl|Y_j\as-E(Y_j\as\mi\cX)\bigr| \si_j(\cX)
>C_7\bigr\}\Bigmi\cX\bigr] ,
\end{eqnarray*}
where $C_6,C_7>0$ are constants, and $C_7$ depends on $\eta$.

Note too that, using (\ref{P28}) to obtain the second relation
below, and (\ref{T26}) to get the last relation, we have
$(\nu\rho h^2)^5\asymp(\nu\delta h^2)^5
%&\asymp(\nu\de)^5 \{(N \nu\de^2)\mo(1-\de)^{-\nu}\}^2\\
\asymp\{(\nu^3 \delta)\half N\mo(1-\delta
)^{-\nu}\}^2
\ra0 .
$
Therefore, (\ref{P22}) holds.
Since $|Y_j\as-E(Y_j\as\mi\cX)|\leq1$, then, if $\si_j(\cX)\leq
C_7$, we have $I\{|Y_j\as-E(Y_j\as\mi\cX)| \si_j(\cX)>C_7\}=0$.
Hence, using~(\ref{P24}) and (\ref{P22}),
\begin{eqnarray*}
S_{N2}(\eta)&\leq&\sum_{j=1}^m
E^0\bigl\{\bigl|Y_j\as-E(Y_j\as\mi\cX)\bigr|^2\Bigmi\cX\bigr\}
\si_j(\cX)^2
I\{\si_j(\cX)>C_7\}\\[-3pt]
&=&(1-\delta)^\nu\{1-(1-\delta)^\nu\} \sum_{j=1}^m
\si_j(\cX)^2
I\{\si_j(\cX)>C_7\}\\[-3pt]
&\leq&(1-\delta)^\nu C_7\mt\sum_{j=1}^m \si_j(\cX)^4
=O_p\{(1-\delta)^\nu m (\nu\delta h^2)^4\}\\[-3pt]
&=&o_p\{m (1-\delta)^\nu(\nu\delta h^2)^2\}
=o_p(1)
\end{eqnarray*}
since $m (1-\delta)^\nu(\nu\delta h^2)^2=C_2$; see (\ref
{P27}). This completes the proof of (\ref{P30}).

\section*{Acknowledgments}

We thank three referees and an Associate Editor for their helpful
comments which led to an improved version of the manuscript.

\begin{supplement}[id=suppA]
\stitle{Additional material}
\slink[doi]{10.1214/11-AOS952SUPP} %[doi,text={...}] - jei reikia
%suskaldyti doi
\sdatatype{.pdf}
\sfilename{aos952\_supp.pdf}
\sdescription{The supplementary article contains a description of
Delaigle and Meister's method, details for bandwidth choice, an
alternative procedure for multivariate setting and unequal groups, and
additional numerical results.}
\end{supplement}

% imsref loaded by lrinkeviciute, 2012-01-25 09:00:14
% imsref loaded by lrinkeviciute, 2012-01-25 09:05:03

\printaddresses


\begin{thebibliography}{21}
% BibTex style file: ims.bst, 2011-05-30
% Default style options (sort=0,type=number).
% Used options (sort=1,type=nameyear).

%b1 ###
%b1 #&#
\bibitem[\protect\citeauthoryear{Bilder and Tebbs}{2009}]{BilTeb09}
\begin{barticle}[mr]
\bauthor{\bsnm{Bilder},~\bfnm{Christopher~R.}\binits{C.~R.}} \AND
  \bauthor{\bsnm{Tebbs},~\bfnm{Joshua~M.}\binits{J.~M.}}
(\byear{2009}).
\btitle{Bias, efficiency, and agreement for group-testing regression models}.
\bjournal{J. Stat. Comput. Simul.}
\bvolume{79}
\bpages{67--80}.
\bid{doi={10.1080/00949650701608990}, issn={0094-9655}, mr={2655675}}
\bptok{imsref}%
\end{barticle}
\endbibitem

%b2 ###
%b2 #&#
\bibitem[\protect\citeauthoryear{Chen and Swallow}{1990}]{CheSwa90}
\begin{barticle}[pbm]
\bauthor{\bsnm{Chen},~\bfnm{C.~L.}\binits{C.~L.}} \AND
  \bauthor{\bsnm{Swallow},~\bfnm{W.~H.}\binits{W.~H.}}
(\byear{1990}).
\btitle{Using group testing to estimate a proportion, and to test the binomial
  model}.
\bjournal{Biometrics}
\bvolume{46}
\bpages{1035--1046}.
\bid{issn={0006-341X}, pmid={2085624}}
\bptok{imsref}%
\end{barticle}
\endbibitem

%b3 ###
%b3 #&#
\bibitem[\protect\citeauthoryear{Chen, Tebbs and Bilder}{2009}]{CheTebBil09}
\begin{barticle}[mr]
\bauthor{\bsnm{Chen},~\bfnm{Peng}\binits{P.}},
  \bauthor{\bsnm{Tebbs},~\bfnm{Joshua~M.}\binits{J.~M.}} \AND
  \bauthor{\bsnm{Bilder},~\bfnm{Christopher~R.}\binits{C.~R.}}
(\byear{2009}).
\btitle{Group testing regression models with fixed and random effects}.
\bjournal{Biometrics}
\bvolume{65}
\bpages{1270--1278}.
\bid{doi={10.1111/j.1541-0420.2008.01183.x}, issn={0006-341X}, mr={2756515}}
\bptok{imsref}%
\end{barticle}\vadjust{\goodbreak}
\endbibitem

%b4 ###
%b4 #&#
\bibitem[\protect\citeauthoryear{Delaigle and Hall}{2011}]{DelHal}
\begin{bmisc}[auto:STB|2012/01/18|07:48:53]
\bauthor{\bsnm{Delaigle},~\bfnm{A.}\binits{A.}} \AND
  \bauthor{\bsnm{Hall},~\bfnm{P.}\binits{P.}}
(\byear{2011}).
\bhowpublished{Supplement to ``Nonparametric regression with homogeneous group
  testing data.'' DOI:\href{http://dx.doi.org/10.1214/11-AOS952SUPP}{10.1214/11-AOS952SUPP}.}
\bptok{imsref}%
\end{bmisc}
\endbibitem

%b5 ###
%b5 #&#
\bibitem[\protect\citeauthoryear{Delaigle and Meister}{2011}]{DelMei11}
\begin{barticle}[auto:STB|2012/01/18|07:48:53]
\bauthor{\bsnm{Delaigle},~\bfnm{A.}\binits{A.}} \AND
  \bauthor{\bsnm{Meister},~\bfnm{A.}\binits{A.}}
(\byear{2011}).
\btitle{Nonparametric regression analysis for group testing data}.
\bjournal{J. Amer. Statist. Assoc.}
\bvolume{106}
\bpages{640--650}.
\bptok{imsref}%
\end{barticle}
\endbibitem

%b6 ###
%b6 #&#
\bibitem[\protect\citeauthoryear{Dorfman}{1943}]{Dor43}
\begin{barticle}[auto:STB|2012/01/18|07:48:53]
\bauthor{\bsnm{Dorfman},~\bfnm{R.}\binits{R.}}
(\byear{1943}).
\btitle{The detection of defective members of large populations}.
\bjournal{Ann. Math. Statist.}
\bvolume{14}
\bpages{436--440}.
\bptok{imsref}%
\end{barticle}
\endbibitem

%b7 ###
%b7 #&#
\bibitem[\protect\citeauthoryear{Fahey, Ourisson and Degnan}{2006}]{FahOurDeg06}
\begin{barticle}[auto:STB|2012/01/18|07:48:53]
\bauthor{\bsnm{Fahey},~\bfnm{J.~W.}\binits{J.~W.}},
  \bauthor{\bsnm{Ourisson},~\bfnm{P.~J.}\binits{P.~J.}} \AND
  \bauthor{\bsnm{Degnan},~\bfnm{F.~H.}\binits{F.~H.}}
(\byear{2006}).
\btitle{Pathogen detection, testing, and control in fresh broccoli sprouts}.
\bjournal{Nutrition J.}
\bvolume{5}
\bpages{13}.
\bptok{imsref}%
\end{barticle}
\endbibitem

%b8 ###
%b8 #&#
\bibitem[\protect\citeauthoryear{Fan}{1993}]{Fan93}
\begin{barticle}[mr]
\bauthor{\bsnm{Fan},~\bfnm{Jianqing}\binits{J.}}
(\byear{1993}).
\btitle{Local linear regression smoothers and their minimax efficiencies}.
\bjournal{Ann. Statist.}
\bvolume{21}
\bpages{196--216}.
\bid{doi={10.1214/aos/1176349022}, issn={0090-5364}, mr={1212173}}
\bptok{imsref}%
\end{barticle}
\endbibitem

%b9 ###
%b9 #&#
\bibitem[\protect\citeauthoryear{Fan and Gijbels}{1996}]{FanGij96}
\begin{bbook}[mr]
\bauthor{\bsnm{Fan},~\bfnm{J.}\binits{J.}} \AND
  \bauthor{\bsnm{Gijbels},~\bfnm{I.}\binits{I.}}
(\byear{1996}).
\btitle{Local Polynomial Modelling and Its Applications}.
\bseries{Monographs on Statistics and Applied Probability}
\bvolume{66}.
\bpublisher{Chapman and Hall}, \baddress{London}.
\bid{mr={1383587}}
\bptok{imsref}%
\end{bbook}
\endbibitem

%b10 ###
%b10 #&#
\bibitem[\protect\citeauthoryear{Fan, Heckman and Wand}{1995}]{FanHecWan95}
\begin{barticle}[mr]
\bauthor{\bsnm{Fan},~\bfnm{Jianqing}\binits{J.}},
  \bauthor{\bsnm{Heckman},~\bfnm{Nancy~E.}\binits{N.~E.}} \AND
  \bauthor{\bsnm{Wand},~\bfnm{M.~P.}\binits{M.~P.}}
(\byear{1995}).
\btitle{Local polynomial kernel regression for generalized linear models and
  quasi-likelihood functions}.
\bjournal{J. Amer. Statist. Assoc.}
\bvolume{90}
\bpages{141--150}.
\bid{issn={0162-1459}, mr={1325121}}
\bptok{imsref}%
\end{barticle}
\endbibitem

%b11 ###
%b11 #&#
\bibitem[\protect\citeauthoryear{Gastwirth and Hammick}{1989}]{GasHam89}
\begin{barticle}[mr]
\bauthor{\bsnm{Gastwirth},~\bfnm{Joseph~L.}\binits{J.~L.}} \AND
  \bauthor{\bsnm{Hammick},~\bfnm{Patricia~A.}\binits{P.~A.}}
(\byear{1989}).
\btitle{Estimation of the prevalence of a rare disease, preserving the
  anonymity of the subjects by group testing: Applications to estimating the
  prevalence of {AIDS} antibodies in blood donors}.
\bjournal{J. Statist. Plann. Inference}
\bvolume{22}
\bpages{15--27}.
\bid{doi={10.1016/0378-3758(89)90061-X}, issn={0378-3758}, mr={0996796}}
\bptok{imsref}%
\end{barticle}
\endbibitem

%b12 ###
%b12 #&#
\bibitem[\protect\citeauthoryear{Gastwirth and Johnson}{1994}]{GasJoh94}
\begin{barticle}[auto:STB|2012/01/18|07:48:53]
\bauthor{\bsnm{Gastwirth},~\bfnm{J.~L.}\binits{J.~L.}} \AND
  \bauthor{\bsnm{Johnson},~\bfnm{W.~O.}\binits{W.~O.}}
(\byear{1994}).
\btitle{Screening with cost-effective quality control: Potential applications
  to HIV and drug testing}.
\bjournal{J. Amer. Statist. Assoc.}
\bvolume{89}
\mbox{\bpages{972--981}}.
\bptok{imsref}%
\end{barticle}
\endbibitem

%b13 ###
%b13 #&#
\bibitem[\protect\citeauthoryear{Hardwick, Page and Stout}{1998}]{HarPagSto98}
\begin{barticle}[mr]
\bauthor{\bsnm{Hardwick},~\bfnm{Janis}\binits{J.}},
  \bauthor{\bsnm{Page},~\bfnm{Connie}\binits{C.}} \AND
  \bauthor{\bsnm{Stout},~\bfnm{Quentin~F.}\binits{Q.~F.}}
(\byear{1998}).
\btitle{Sequentially deciding between two experiments for estimating a common
  success probability}.
\bjournal{J. Amer. Statist. Assoc.}
\bvolume{93}
\bpages{1502--1511}.
\bid{issn={0162-1459}, mr={1666644}}
\bptok{imsref}%
\end{barticle}
\endbibitem

%b14 ###
%b14 #&#
%  \bauthor{\bsnm{Swallow},~\bfnm{William~H.}\binits{W.~H.}}
%(\byear{2000}).
%  or quantitative covariables}.

%b15 ###
%b15 #&#
\bibitem[\protect\citeauthoryear{Lennon}{2007}]{Len07}
\begin{barticle}[pbm]
\bauthor{\bsnm{Lennon},~\bfnm{Jay~T.}\binits{J.~T.}}
(\byear{2007}).
\btitle{Diversity and metabolism of marine bacteria cultivated on dissolved
  DNA}.
\bjournal{Applied and Environmental Microbiology}
\bvolume{73}
\bpages{2799--2805}.
\bid{doi={10.1128/AEM.02674-06}, issn={0099-2240}, pii={AEM.02674-06},
  pmcid={1892854}, pmid={17337557}}
\bptok{imsref}%
\end{barticle}
\endbibitem

%b16 ###
%b16 #&#
\bibitem[\protect\citeauthoryear{Nagi and Raggi}{1972}]{NagRag}
\begin{barticle}[pbm]
\bauthor{\bsnm{Nagi},~\bfnm{M.~S.}\binits{M.~S.}} \AND
  \bauthor{\bsnm{Raggi},~\bfnm{L.~G.}\binits{L.~G.}}
(\byear{1972}).
\btitle{Importance to ``airsac'' disease of water supplies contaminated with
  pathogenic Escherichia coli}.
\bjournal{Avian Diseases}
\bvolume{16}
\bpages{718--723}.
\bid{issn={0005-2086}, pmid={4562575}}
\bptok{imsref}%
\end{barticle}
\endbibitem

%b17 ###
%b17 #&#
\bibitem[\protect\citeauthoryear{Ruppert, Wand and Carroll}{2003}]{RupWanCar03}
\begin{bbook}[mr]
\bauthor{\bsnm{Ruppert},~\bfnm{David}\binits{D.}},
  \bauthor{\bsnm{Wand},~\bfnm{M.~P.}\binits{M.~P.}} \AND
  \bauthor{\bsnm{Carroll},~\bfnm{R.~J.}\binits{R.~J.}}
(\byear{2003}).
\btitle{Semiparametric Regression}.
\bseries{Cambridge Series in Statistical and Probabilistic Mathematics}
\bvolume{12}.
\bpublisher{Cambridge Univ. Press}, \baddress{Cambridge}.
\bid{doi={10.1017/CBO9780511755453}, mr={1998720}}
\bptok{imsref}%
\end{bbook}
\endbibitem

\bibitem[\protect\citeauthoryear{Thorburn et~al.}{2001}]{Thoetal01}
\begin{barticle}[pbm]
\bauthor{\bsnm{Thorburn},~\bfnm{D.}\binits{D.}},
  \bauthor{\bsnm{Dundas},~\bfnm{D.}\binits{D.}},
  \bauthor{\bsnm{McCruden},~\bfnm{E.}\binits{E.}},
  \bauthor{\bsnm{Cameron},~\bfnm{S.}\binits{S.}},
  \bauthor{\bsnm{Goldberg},~\bfnm{D.}\binits{D.}},
  \bauthor{\bsnm{Symington},~\bfnm{I.}\binits{I.}},
  \bauthor{\bsnm{Kirk},~\bfnm{A.}\binits{A.}} \AND
  \bauthor{\bsnm{Mills},~\bfnm{P.}\binits{P.}}
(\byear{2001}).
\btitle{A study of hepatitis C prevalence in healthcare workers in the
west of Scotland}.
\bjournal{Gut}
\bvolume{48}
\bpages{116--120}.
\bid{issn={0017-5749}, pmcid={1728181}, pmid={11115832}}
\bptok{imsref}%
\end{barticle}
\endbibitem

%b18 ###
%b18 #&#
\bibitem[\protect\citeauthoryear{Vansteelandt, Goetghebeur and Verstraeten}{2000}]{VanGoeVer00}
\begin{barticle}[mr]
\bauthor{\bsnm{Vansteelandt},~\bfnm{S.}\binits{S.}},
  \bauthor{\bsnm{Goetghebeur},~\bfnm{E.}\binits{E.}} \AND
  \bauthor{\bsnm{Verstraeten},~\bfnm{T.}\binits{T.}}
(\byear{2000}).
\btitle{Regression models for disease prevalance with diagnostic tests on pools
  of serum samples}.
\bjournal{Biometrics}
\bvolume{56}
\bpages{1126--1133}.
\bid{doi={10.1111/j.0006-341X.2000.01126.x}, issn={0006-341X}, mr={1806746}}
\bptok{imsref}%
\end{barticle}
\endbibitem

%b19 ###
%b19 #&#
\bibitem[\protect\citeauthoryear{Wahed et~al.}{2006}]{Wahetal06}
\begin{barticle}[auto:STB|2012/01/18|07:48:53]
\bauthor{\bsnm{Wahed},~\bfnm{M.~A.}\binits{M.~A.}},
  \bauthor{\bsnm{Chowdhury},~\bfnm{D.}\binits{D.}},
  \bauthor{\bsnm{Nermell},~\bfnm{B.}\binits{B.}},
  \bauthor{\bsnm{Khan},~\bfnm{S.~I.}\binits{S.~I.}},
  \bauthor{\bsnm{Ilias},~\bfnm{M.}\binits{M.}},
  \bauthor{\bsnm{Rahman},~\bfnm{M.}\binits{M.}},
  \bauthor{\bsnm{Persson},~\bfnm{L.~A.}\binits{L.~A.}} \AND
  \bauthor{\bsnm{Vahter},~\bfnm{M.}\binits{M.}}
(\byear{2006}).
\btitle{A modified routine analysis of arsenic content in drinking-water in
  Bangladesh by hydride generation-atomic absorption spectrophotometry}.
\bjournal{J. Health, Population and Nutrition}
\bvolume{24}
\bpages{36--41}.
\bptok{imsref}%
\end{barticle}
\endbibitem

%b20 ###
%b20 #&#
\bibitem[\protect\citeauthoryear{Xie}{2001}]{Xie01}
\begin{barticle}[auto:STB|2012/01/18|07:48:53]
\bauthor{\bsnm{Xie},~\bfnm{M.}\binits{M.}}
(\byear{2001}).
\btitle{Regression analysis of group testing samples}.
\bjournal{Stat. Med.}
\bvolume{20}
\bpages{1957--1969}.
\bptok{imsref}%
\end{barticle}
\endbibitem

\end{thebibliography}
\end{document}